\documentclass[11pt]{report}
\usepackage[centertags]{amsmath}
\usepackage{amsfonts}
\usepackage{amssymb}
\usepackage{amsthm}
\usepackage{amscd}
\usepackage[all]{xy}
\usepackage{url}

\usepackage{graphicx}

\usepackage{amssymb}
\usepackage{amsfonts}
\usepackage{enumerate}
\usepackage{amsthm}

\newcommand{\suml}{\sum\limits}
\newcommand{\vare}{\varepsilon}
\newcommand{\bbr}{\mathbb{R}}
\newcommand{\bbn}{\mathbb{N}}
\newcommand{\bbz}{\mathbb{Z}}
\newcommand{\bbf}{\mathbb{F}}
\newcommand{\bbh}{\mathbb{H}}
\newcommand{\bbq}{\mathbb{Q}}
\newcommand{\calb}{\mathcal{B}}
\newcommand{\calp}{\mathcal{P}}
\newcommand{\ga}{\mathfrak{A}}
\def\thalf{{\textstyle{\frac 12}}}
\newtheorem{lemma}{Lemma}[chapter]
\newtheorem{lem}[lemma]{Lemma}
\newtheorem{theorem}[lemma]{Theorem}
\newtheorem{rem}[lemma]{Remark}
\theoremstyle{definition}
\newtheorem{cor}[lemma]{Corollary}
\newtheorem{conj}[lemma]{Conjecture}
\newtheorem{definition}[lemma]{Definition}
\newtheorem{prop}[lemma]{Proposition}
\newtheorem{exa}[lemma]{Example}

\newtheorem{open}[lemma]{Open Problem}

\def\be*{\begin{equation*}}
\def\ee*{\end{equation*}}
\def\call{\mathcal{L}}
\def\caln{\mathcal{N}}
\def\calh{\mathcal{H}}
\def\la{\lambda}
\def\a{\alpha}
\def\intl{\int\limits}
\def\tle{\triangleleft}
\def\ol{\overline}
\def\Ga{\Gamma}

\newcommand{\bbc}{\mathbb{C}}

\newcommand{\cO}{{\mathcal O}}

\newcommand{\bN}{{\mathbb N}}

\newcommand{\bR}{{\mathbb R}}
\newcommand{\bC}{{\mathbb C}}
\newcommand{\bZ}{{\mathbb Z}}
\newcommand{\bQ}{{\mathbb Q}}

\newcommand{\bP}{{\mathbb P}}

\begin{document}
\begin{titlepage}

\begin{center}


{ \huge \bfseries Expander Graphs\\
\ \\
\Large in Pure and Applied Mathematics}\\[1.5cm]
\centerline{\large{\bfseries  Alexander Lubotzky}}

\smallskip

\centerline{\bf Einstein Institute of Mathematics}
\centerline{\bf  Hebrew University}
\centerline{\bf  Jerusalem 91904}
\centerline{\bf  ISRAEL}
\centerline{\bf alexlub@math.huji.ac.il}

\vfill
\end{center}
\underline {Abstract:} \\
Expander graphs are highly connected sparse finite graphs. They play an important role in computer science as basic building blocks for network constructions, error correcting codes, algorithms and more. In recent years they have started to play an increasing role also in pure mathematics: number theory, group theory, geometry and more. This expository article describes their constructions and various applications in pure and applied mathematics.
\vfill
\begin{center}
\noindent {\large This paper is based on notes prepared for the Colloquium Lectures at the\\  Joint Annual Meeting of the American Mathematical Society (AMS)\\ and the Mathematical Association of America (MAA).\\
New Orleans, January 6-9, 2011}.\\ \ \\ The author is grateful to the AMS for the opportunity to present this material for a wide audience. He has benefited by responses and remarks \\ which followed his lectures.\\ \ \\
\underline {2010 Math Subject Classification:} 01-02, 05C99
\end{center}

\end{titlepage}

\tableofcontents

\vfill

\noindent {\bf Acknowledgments.} The author is indebted to Peter Sarnak for many years of fruitful collaboration and friendship.  Much of what is presented here was inspired by him and in particular in Chapter 4 we have made extensive use of material from his web-site.  We are grateful to E. Kowalski, N. Linial and C. Meiri for helpful comments on an earlier draft.
Thanks are  also due to the ERC and ISF for partial support.

\bigskip

\newpage

 \chapter*{Introduction}
Expander graphs are highly connected sparse finite graphs which play a basic role in various areas of computer science. A huge amount of research has been devoted to them in the computer science literature in the last four decades. (An excellent survey of these directions is \cite{HLW}). But they also attracted the attention of mathematicians: Their existence follows easily by random considerations ($\grave{a}$ $la$ Erd\H{o}s) but explicit constructions, which are very desirable for applications, are much more difficult. Various deep mathematical theories have been used to give explicit constructions, e.g. Kazhdan property ($T$) from representation theory of semi-simple Lie groups and their discrete subgroups, the Ramanujan Conjecture (proved by Deligne) from the theory of automorphic forms and more. All these led to fascinating connections between pure mathematics and computer science and between pure mathematicians and computer scientists. For the first three decades -- till approximately ten years ago -- essentially all these connections went in one direction: methods of pure mathematics have been used to solve some problems arising from computer science (these are summarized in \cite{L1} for example).

Something different is emerging in the last decade: Computer science pays its debt to pure mathematics! The notion of expander graphs is starting to play a significant role in more and more areas of pure mathematics.

The goal of these notes is to describe expander graphs and their applications in pure and applied mathematics. Rather than competing with the award winning manuscripts \cite{L1} and \cite{HLW} (the Ferran Sunyer i Balaguer prize and the Levi Conant prize, respectively), we will place  emphasis on the new directions: applications of expanders in pure mathematics. We will try to avoid repeating topics from \cite{L1} and \cite{HLW}, though some intersection is unavoidable, especially in the first chapters. The reader is strongly encouraged to consult these manuscripts  for more background, as well as \cite{LZu}.

The notes are organized as follows. We start with basic definitions of expander graphs, their properties, their eigenvalues and random walks on them. In the second chapter we will give various examples, mainly of Cayley graphs, which are expanders. The reader should not be misled by this chapter's modest title, ``examples". Some of the most remarkable developments in recent years are described there, e.g., the fact that all non-abelian finite simple groups are expanders in a uniform way and the result that congruence quotients of linear groups form a family of expanders. The last result is the crucial ingredient in some of the applications to number theory and to group theory.

Chapter 3 deals with applications to computing. Many are described in \cite{HLW}, so we chose to give a theoretical application to the product replacement algorithm and one for error correcting codes.

Chapter 4 deals with applications to number theory. There are several of these, but we will mainly describe a new direction of research for which the use of expanders is a dominant factor: ``the affine sieve". This method enables to study primes and almost primes in orbits of groups acting on $\bZ^n$. This is a far-reaching extension to a non-commutative world of Dirichlet's theorem about primes in arithmetic progression. This direction of research arose in response to the dramatic developments concerning Cayley graphs being expanders. It also sheds new light on classical subjects like Apollonian circle packing and more.

The affine sieve method can be also modified to give ``group sieve" which is a method to study various group theoretical properties of ``generic" elements in finitely generated groups. This method gives some new results about linear groups and the mapping class groups and these are described in chapter 5.

Chapter 6 is devoted to applications to geometry. Most of the applications are for hyperbolic manifolds with some special attention to hyperbolic 3-manifolds.

In chapter 7, we collected brief remarks on several topics which should fit into these notes but for various reasons were left out.

We hope that the current notes will provide a panoramic view of the broad scope of mathematics which is connected with expander graphs. They have truly expanded into many different areas of mathematics!

\chapter{Expander graphs}

Expander graphs are highly connected sparse graphs.  This property can be viewed from several different angles: eigenvalues, random walks, representation theory (if the graph is a Cayley graph), geometry and more.  In this chapter we briefly review these aspects (sending the reader to [L1] for a more comprehensive description).  This leads to the highly relevant property $(\tau)$ described in \S 1.6, which will play a very important role in the chapters to come.

\section{The basic definition}
 Let $X$ be a finite graph on a set $V$ of $n$ vertices  and $A=A_X$ its adjacency matrix, i.e. $A$ is an $n\times n$ matrix where $A_{i, j}$ is the number of edges between vertex $i$ and vertex $j$.  So usually $A_{i, j} = 0 $ or $1$, but we also allow multiple edges $(A_{ij} > 1)$ or loops $(A_{ii}> 0)$.
 The graph $X$ is $k$-{\it regular} if the valency of every vertex is $k$, i.e. for every $i$, $\sum^n_{j=1} A_{ij} = k$.

 \begin{definition}\label{def1.1} For $0 < \vare \in\bbr, \; X$ is $\pmb\vare${\bf-expander} if for every subset $Y $ of $V$
 with
 \begin{equation*}|Y| \le \thalf |V| = \frac n2, \; \; \; \; |\partial Y|\ge \vare |Y|\end{equation*} where $\partial Y$ is the boundary of $Y$,  i.e., the set of vertices in $V$ which are connected to (some vertices of) $Y$ but are not in $Y$.

The largest $\vare$ for which $X$ is $\vare$-expander will be denoted $\vare (X)$.
\end{definition}

One easily sees that $X$ is connected. So being $\vare$-expander for ``large" $\vare$ (well, it is clear that $\vare$ cannot be larger than 1) means that $X$ is ``very much connected".

In most applications (real world applications as well as pure mathematical applications) what one wants is to find regular graphs with large $n$ (say $n\to\infty)$, fixed $k$ (as small as possible) and a fixed $\vare $ (as large as possible). A family of $k$-regular graphs will be called a \emph{family of expanders} or an \emph{expanding family} if all of them are $\vare$-expanders for the same $\vare > 0$.

The first to define expander graphs was Pinsker [Pin] in 1973 who also coined the name.  Recently, it has been noticed by Larry Guth that slightly earlier Barzdin and Kolmogorov [BK] discussed a property of graphs which is equivalent to expanders.  While Pinsker defined and studied expanders for their use in computer science (error correction codes, communication networks and algorithms) Barzdin and Kolmogorov's motivation was very different: they studied the network of nerve cells of the human brain.  This brought them to the question on realizing various networks in $\bbr^3$ and this way to graphs in which any two subsets have a large number of edges between them, a property which characterizes expanders (see [HLW, \S 2.4], and the historical notes in [GrGu]).

\smallskip

When $k$ is fixed, bounding $\vare$ - the expansion constant is closely related to the \emph{isoperimetric constant} $h(X)$ called also - the \emph{Cheeger constant}.

\begin{definition}\label{def1.2}
For $X$ as above, let
\be* h(X):= \min\limits_{V = Y_1\dot\bigcup
Y_2} \; \frac{|E(Y_1, Y_2)|}{\min((|Y_1|, |Y_2|)}
\ee*
where the minimum runs over all the ways to write $V$ as a disjoint union of two subsets $Y_1$ and $Y_2$ and $E(Y_1, Y_2)$ is the set of edges between $Y_1$ and $Y_2$.
\end{definition}

The following is a straight-forward corollary  of the definitions:

\begin{prop}\label{prop1.1}
\be* \frac{h(X)}{k} \le \vare (X) \le h (X)\ee*
\end{prop}

\section{Eigenvalues, random walks and Ramanujan graphs}

Let $X$ be as before, a $k$-regular graph on $n$ vertices.  As $X$ is undirected, $A=A_X$ is a symmetric matrix with real eigenvalues.  One can think of $A$ as a linear operator on $L^2(X)$ - the real (or complex) functions on $V$, where for $i \in V$ and $f \in L^2(X)$,
\be* (Af) (i) = \suml^n_{j=1} A_{ij} f(j)\ee*
i.e., summing $f$ on the neighbors.
An easy argument shows that $k$ is the largest eigenvalue of $A$ (corresponding to the constant functions) and all the eigenvalues $\la_0 = k \ge \la_1 \ge \dots \ge \la_{n-1}$ of $A_X$ lie in the interval $[-k,k]$.  Some well known easy properties:
\begin{prop}\label{prop1.2} (a) $X$ is connected iff $\la_1 < \la_0 = k$

(b) $X$ is bipartite iff $\la_{n-1} = - k$.
\end{prop}

So various combinatorial/geometric properties of $X$ can be recovered from its spectrum (eigenvalues).  So is the expansion property:

\begin{prop}\label{prop1.3}
\be* \frac{k-\la_1}{2} \le h (X) \le \sqrt{(k+\lambda_1)(k-\la_1)}\ee*

\end{prop}

So for $k$-regular graphs (fixed $k$) being $\vare$-expanders is equivalent to a spectral gap $\la_1 < k - \vare'$.

 This proposition is now well known and is attributed to various authors
(see [L1] for more details and references about this result and the other ones in this chapter).

For the original definition of  expander one needs to bound the largest eigenvalue $\lambda_1 = \lambda_1(X)$  of $X$ which is smaller than $k$.  But for various other applications what is most relevant is:
\be*
\la(X): = \max\{ |\la| \big| \lambda \text{\ an\ eigenvalue\ of\ $A$ and\ } |\la|\neq k\}
\ee*
i.e. the largest eigenvalue in absolute value other than $\pm k$. This is not a crucial difference  but some care is needed as some authors define expanders by a bound on $\la(X)$.

For many applications the bound on the eigenvalues is even more relevant than the original definition since the eigenvalues control the random walk on $X$.  Namely, assume $\mu \in L^2(X)$ is a probability measure on $V$, i.e. $0 \le \mu (i) \le 1$ for every $i\in V$ and $\sum^n_{i = 1} \mu(i) = 1$.  Then if a ``little person" is in vertex $i$ at step $t$ with probability $\mu(i)$ and then he walks a step over a randomly chosen edge coming out of $i$, then at step $t + 1$ he will be in $\frac 1k A\mu$.  In other words the matrix $\triangle = \triangle_X = \frac 1k A_X$ is the bi-stochastic transition  matrix of the Markov chain that is the random walk on  $X$.  If $X$ is connected and not bi-partite then the random walk converges to the uniform distribution $u$, i.e. $u(i) = \frac 1n \; \text{for\ every\ } i \in V$. The rate of convergence depends on
\be* \la(X):= \max \{ | \la | \,\big| \, |\la| \neq k,\quad \la \text{\ an\ eigenvalue\ of\ } A\}\ee*

More precisely (cf. \cite[Theorem 3.3]{HLW}),
\begin{prop}\label{pr1.6}
Let $X$ be a non-bi-partite $k$-regular graph with adjacency matrix $A$ and normalized one $\Delta = \frac 1k A$. Then for any distribution $\mu$ on the vertices of $X$ and any $1 \le t \in\bbn$,
\be* \| \Delta^t \mu - u\|_{L_2} \le \big(\frac{\la(X)}{k}\big)^t\ee*
when $u$ is the uniform distribution.
\end{prop}

The non-bi-partite issue is not crucial and can be avoided by considering the ``lazy random walk" - cf. \cite{LP}.  One can also get similar types of bounds with the $L^1$-norm which is sometimes more relevant  - see \cite[\S 3.1]{HLW}.

 There is a limit to what one can expect when trying to bound $\la(X)$.  This is given by the Alon-Boppana result:

\begin{prop}\label{prop1.4} Let $X_{n, k}$ be an infinite family of $k$-regular connected graphs on $n$ vertices where $k$ is fixed and $n\to \infty$.

Then $ \la(X_{n, k}) \ge 2 \sqrt{k-1} - o(1)$.
\end{prop}

This suggests the following definition:

\begin{definition}\label{def1.3}  A $k$-regular finite graph $X$ is called  a \emph{Ramanujan graph} if $\la(X) \le 2 \sqrt{k-1}$.
\end{definition}

So, Ramanujan graphs are, in some sense, optimal expanders.
The most general known result gives for every $k$ of the form $p^\a +1$, where $p$ is a prime and $\a \in\bbn$, an infinite family of $k$-regular Ramanujan graphs (\cite{Mo}, \cite{LSV2}).  We mention in passing that for every $k$ which is not of this form, it is not known if such an infinite family exists.  The first open case is $k = 7$.

\section{Cayley graphs and representation theory}

A particularly nice way to construct graphs which are very symmetric is via Cayley graphs.  Recall that if $G$ is group and $\Sigma$ a symmetric subset of $G$ (i.e., $ s \in \Sigma$ iff $s^{-1} \in \Sigma$), the Cayley graph $Cay (G; \Sigma)$ of $G$ w.r.t. $\Sigma $ is the graph whose vertex set is $G$ and $a\in G$ is connected to $\{ s a | s \in \Sigma\}$.  This is a $k$-regular graph with $k=|\Sigma|$.  It is connected iff $\Sigma$ generates $G$.

The expansion properties of $Cay (G; \Sigma)$ can be reformulated in representation theoretic terms.  To this end, let us define:
\begin{definition}\label{def1.4} Let $G$ be a group and $\Sigma$ a subset of $G$. We say that $\vare' > 0$ is a \emph{Kazhdan constant} of $G$ w.r.t. $\Sigma$ if for every unitary representation $\rho: G\to U(H)$, where $H$ is a Hilbert space and $U(H)$ the group of unitary operators, without a non-zero fixed vector, and for every $0\neq v \in H$, there exists $s \in \Sigma \ \text{\ such\ that\ } \| \rho(s) v-v\| \ge \vare' \| v \|$.
\end{definition}

One can see that in this case $\Sigma$ generates $G$.

\begin{definition}\label{def1.5} A discrete group $\Gamma$ is said to have \emph{Kazhdan property $(T)$}, if it has some finite set of generators $\Sigma$ with Kazhdan constant $\vare' > 0$.

\end{definition}

One can prove that if this happens for one $\Sigma$ it is so for any set of generators, with possibly different $\vare'$.

The Kazhdan constant is another way to express the expansion of Cayley graphs.

\begin{prop}\label{15}  (i) For every $0 < \vare' \in\bbr$, there exists $\vare = f_1 (\vare')> 0 $ s.t. if $G$ is a finite group with a symmetric set of generators $\Sigma$ and Kazhdan constant $\vare'$, then $\vare(Cay (G;\Sigma)) \ge \vare$, i.e., $Cay (G; \Sigma)$ is an $\vare$-expander.

(ii)  For every $k \in \bbn$ and every $0 < \vare \in \bbr$, there exists $\vare' = f_2 (k, \vare)$, such that if $G$ is a finite group with a symmetric set  of $k$ generators $\Sigma$ with $\vare(Cay (G;\Sigma)) \ge \vare$, $\vare' = f_2 (k, \vare)$ is a Kazhdan constant for $G$ w.r.t. $\Sigma$.
\end{prop}

So, at least as long as $k$ is fixed the expansion constant and the Kazhdan constant are closely related.  Assume now that $\Gamma$ is an infinite group generated by a finite symmetric set $\Sigma$ and assume $\call = \{ N_i\}_{i\in I}$ is an infinite collection of finite index normal subgroups of $\Gamma$.  We can deduce from Proposition \ref{15}:

\begin{prop}[\cite{M1}]\label{16} If $\Gamma$ has Kazhdan property $(T)$, i.e., there exists $\vare' > 0$ which is a Kazhdan constant of $\Gamma$ w.r.t. $\Sigma$, then all the finite quotients $Cay (\Gamma/N_i; \Sigma)$, $i \in I$, are $\vare$-expander where $\vare > 0$ depends only on $\vare'$.

\medskip

The fact that there are groups with property $(T)$ is a non-trivial result due originally to Kazhdan (see \S 2.2 below for more).  For example $\Gamma = SL_3 (\bbz)$, the integral $3\times 3$ matrices of determinant 1 is such a group and so the Cayley graphs of its quotients $SL_3(\bbz/m\bbz)$, $m\in\bbn$, form a family of $\vare$-expanders w.r.t. a fixed set of generators coming from $SL_3 (\bbz)$, e.g. $\{ A^{\pm 1}, B^{\pm 1}\}$ where
\be* A= \Bigl(\begin{smallmatrix} 1 &1 &0\\
0 &1 &0\\
0 &0 &1\end{smallmatrix}\Bigr) \quad \text{and}\quad  B = \Bigl(\begin{smallmatrix} 0 &1 &0\\
0 &0 &1\\
1 &0 &0\end{smallmatrix}\Bigr).
\ee*
\end{prop}

We will come back to many more examples of this kind in Chapter 2. Here we only observe that in order to deduce the conclusion of Proposition \ref{16}, we need a weaker property than $(T)$, the so called $(\tau)$.  Due to  its importance we will give the definition repeating the notations again:

\begin{definition}\label{def1.6}  Let $\Gamma$ be a group, with a collection $\call = \{ \caln_i\}_{i\in I}$ of finite index normal subgroups.  We say that $\Gamma$ has {\bf property}
{\mbox{(\boldmath{$\tau$}})} w.r.t. $\call$, if there exists a symmetric subset $\Sigma$ of $\Gamma$ and an $0 < \vare' \in\bbr$ such that for every finite quotient $G_i = \Gamma/N_i$, $ i \in I, \; \vare'$ is a Kazhdan constant for $G_i$ with respect to $\Sigma$ (or more precisely w.r.t. $\Sigma N_i/ N_i)$.  An equivalent way to say it is that for every unitary representation $\rho:\Gamma \to U(\calh)$,    with ${ Ker\rho \supset N_i}$ for some ${i}$,
without a non-zero fixed vector, and every $0 \neq v \in H$, there exists $s\in\Sigma$ such that $\| \rho (s) v - v\| > \vare' \| v\|$.
If $\call$ is the family of \emph{all}  finite index normal subgroups of $\Gamma$, we simply say that $\Gamma$ has property $(\tau)$.
\end{definition}

The equivalence of Proposition \ref{15} shows:
\begin{prop}\label{1.7} The group  $\Gamma$ has property $(\tau)$ w.r.t. $\call=\{ N_i\}_{i\in I}$ if and  only if there exists a symmetric subset $\Sigma$ of $\Gamma$ and $\vare > 0$, such that all the Cayley graphs $Cay(\Gamma/N_i; \Sigma)$ are $\vare$-expanders.
\end{prop}

There exist groups (e.g. $\Gamma = SL_2(\bbz[\frac 1p])$) which have $(\tau)$ but not $(T)$  while $\Gamma = SL_2(\bbz)$ has  neither $(T)$ nor $(\tau)$, but it has $(\tau)$ w.r.t. the family
\be* \call = \{ \Gamma (m) = Ker \big(SL_2(\bbz) \to SL_2(\bbz/m\bbz)\big)\}_{m\in\bbn}\ee*
the family of congruence subgroups - see \S 2.4 below.

Some of the recent breakthroughs  are far-reaching extensions of this last fact; extensions which have some remarkable applications.

\section{Expanders and Riemannian manifolds}

Let $M$ be an $n$-dimensional connected closed Riemannian manifold (i.e. compact with no boundary; much of theory can be extended to a more general setting but for simplicity of the exposition we will stick to the closed case).  Let $\Delta = - div (grad)$ be the Laplacian operator of $L^2(M)$.  Its eigenvalues $0 = \lambda_0 (M) < \lambda_1 (M) \le \lambda_2 (M) \le \cdots $ form a discrete subset (with multiplicities) of $\bbr_+$, called the \emph{spectrum} of $M$.

The spectrum of $M$ is very much related to the geometry of $M$ and these relations are the subject of \emph{spectral geometry}. A more intuitive description of $\Delta$ is given by the formula:
\be*
(\Delta f) (p) = \lim\limits_{r\to 0} \frac{2n}{r^2} \Big( \frac{\int_{S^r} f}{vol (S_r)} - f(p)\Big) \ee*
where $n = \dim M, \; p \in M, \; f \in L^2(M)$ and $S_r$ is the sphere of radius $r$ around $p$.  This description  is similar to the combinatorial Laplacian  as an averaging operator.

We will mainly be interested in $\lambda_1 (M)$ which can be described directly without a reference to $\Delta$.

\begin{prop}\label{prop1.8} \be* \lambda_1(M) = \inf \large\{ \frac{\int_M\| df\|^2}{\int_M |f|^2} \big| f \in C^\infty (M), \intl_M f=0\big\}\ee*
\end{prop}

Another important geometric invariant of $M$, whose connection with expanders is even more evident is the Cheeger constant:
\begin{definition}\label{def1.7} The \emph{Cheeger constant} $h(M)$ is
\be* h(M) = \underset{E}{\inf} \, \frac{\mu(E)}{\min(\nu(A), \nu(B))}\ee*
where $E$ runs over all the compact $(n-1)$-dimensional submanifolds of $M$ which divide $M$ into disjoint submanifolds $A$ and $B$.  Here $\mu(E)$ is the ``area" of $E$ and $\nu$ the volume form of $M$.
\end{definition}

Just as for graphs, $h(M)$ is closely related to $\lambda_1 (M)$.  In fact, historically the relation between expansion/Cheeger constant and $\lambda_1$ was discovered first for manifolds and only later on for graphs - see [L1] for historical notes.

\begin{theorem}[Cheeger's inequality]\label{11} $\lambda_1(M)\ge \frac{h^2(M)}{4}$.
\end{theorem}

Buser  proved a converse to this inequality, which depends on the Ricci curvature $R(M)$. We will not bother defining it here - but rather send  the reader to \cite{Bu}, \cite{L1} and the references therein. But let's quote:
\begin{theorem}\label{12}  If $R(M) \ge - (n-1) a^2$ for some $a\ge 0$ where $n=\dim M$, then $\lambda_1(M) \le 2 a (n-1) h (M) + 10 h^2(M)$.
\end{theorem}
 What is important for us is that in case of bounded Ricci curvature, which will hold in all our considerations, $\la_1 (M)$ is also bounded above by a function of $h(M)$.

The more precise connection between these notions and expander graphs will be given in Theorem \ref{14} below.  Let us point out here the basic intuition:

Let $\tilde M$ be the universal cover of $M, \; \Gamma = \pi_1(M)$ the fundamental group of $M$ and $F$ a fundamental domain for the action of $\Gamma$ on $\tilde M$, i.e., $F$ is an open subset of $\tilde M$, whose closure $\bar F$ is compact and such that $\Gamma  \bar F = \tilde M$ and for every $1 \neq \gamma \in \Gamma, \; \; \gamma (\bar F) \cap \bar F \subseteq \bar F\smallsetminus F$.
Standard covering theory shows that the finite set $\Sigma = \{ \gamma\in \Gamma \big| \gamma\bar F\cap \bar F$ is of codimension 1$\}$ is a symmetric set of generators for $\Gamma$.  One can visualize $Cay(\Gamma; \Sigma)$ in the following way: Fix $x_0 \in F$ and put a  vertex at the interior point $\gamma x_0$ of the tesselate $\gamma F$ of $F$ (naturally this vertex will represent $\gamma$; note that $\gamma$ is unique for the given tesselate).  Now, draw an edge between $\gamma_1 x_0$ and $\gamma_2 x_0$ if $\overline{\gamma_1 F} \cap \overline{\gamma_2 (F)}$ is of codimension 1.  One can easily check that what we get is exactly a ``drawing" of $Cay (\Gamma; \Sigma)$ on $\tilde M$.

Moreover, if $\Gamma_1$ is a normal subgroup of $\Gamma$ of finite index, then the ``projection" of the above graph to $M_1 = \tilde M/\Gamma_1$ is exactly the Cayley graph $Cay (\Gamma/\Gamma_1; \Sigma)$.  We therefore get that  the combinatorial graphs $Cay (\Gamma/\Gamma_1; \Sigma)$ when $\Gamma_1$ runs over the finite index normal subgroups of $\Gamma$, are ``approximations" of the finite sheeted normal covers of $M$.

This enables us to relate the expansion properties of these Cayley graphs to the asymptotic of $h(M_1)$ and similarly with $\lambda_1$ of the graphs and of the manifolds - see Theorem \ref{14} below.

\section{Expanders and measure theory}

Let $G$ be a compact group.  A \emph{mean} $m$ on $G$ is a linear functional $m:L^\infty (G)\to\bbr$ satisfying (i) \ $m(f)\ge 0$ if $f \ge 0$; (ii) $ m (\chi_G) = 1 $ where $\chi_G$ is the constant function $1$ on $G$. We say that it is a $\mathbf{G}${\it-invariant mean} if it also satisfies: (iii) \ $m (g.f) = m(f)$ for every $g \in G$ and $f \in L^\infty (G)$, where $g.f (x) = f(g^{-1} x)$, i.e., $G$-left invariant.

An example of an invariant mean is the Haar measure of $G$ which is also  a countable additive on subsets of $G$.  In general, it is possible for  $G$ to  have invariant means different than the (unique) Haar measure.  For example $G=S^1$-the circle has such means. But:
\begin{theorem}[\cite{Sh1}]\label{13}  Let $\Gamma$ be a finitely generated group and $\call = \{ N_i\}_{i \in \bbn}$ a decreasing chain of finite index normal subgroups of $\Gamma$. Let\\
 $G = \lim\limits_{\underset{i\in\bbn}{\longleftarrow}} \Gamma/N_i$ be the profinite completion of $\Gamma$ w.r.t. $\call$.  Then the following are equivalent:
\begin{enumerate}[(i)]
\item $\Gamma$ has $(\tau)$ w.r.t. $\call$

\item The Haar measure on $G$ is the only $\Gamma$-invariant mean on $G$.
   \end{enumerate}
\end{theorem}
So, beside all the previous combinatorial, geometric and representation theoretic characterizations of property $(\tau)$, we also have a measure theoretic one.

\section{A summary of Property ($\tau$)}

Let $\Gamma$ be a finitely generated group generated by a finite symmetric set $\Sigma$.  Let $\call = \{ N_i\}_{i \in \bbn}$ be a decreasing chain of finite index normal subgroups of $\Gamma$.

\begin{theorem}\label{14}
The following conditions are equivalent:
\begin{enumerate}[(i)]
\item $\Gamma$ has property $(\tau)$ w.r.t. $\call$, i.e., there exists $\vare_1 > 0$ s.t. if $\rho: \Gamma \to U(H)$ is a unitary representation of $\Gamma$ on a Hilbert space $H$ without non-zero $\rho(\Gamma)$ fixed points and such that $Ker (\rho) \ge N_i$ for some $i$, then for every $0\neq v \in H$, there exists $s \in\Sigma$ with $\| \rho (s) v-v\| \ge \vare_1 \| v \|$.

    \item There exists $\vare_2 > 0$ such that all the Cayley graphs $Cay (\Gamma/N_i; \Sigma)$ are $\vare_2$-expanders.

        \item There exists $\vare_3 > 0 $ such that $h (Cay (\Gamma/N_i); \Sigma)) \ge \vare_3$ (see Definition \ref{def1.2}).

            \item There exists $\vare_4 > 0$ such that $\lambda_1 (Cay(\Gamma/N_i; \Sigma)) \le k - \vare_4$ for every $i \in\bbn$, where $k = |\Sigma |$.

\item The Haar measure of $\hat\Gamma_{\call} = \underset{\longleftarrow}{\lim} \; \Gamma/N_i$ is the only $\Gamma$-invariant mean on $L^\infty(\hat\Gamma_\call)$.


\medskip

\noindent If in addition $\Gamma = \pi_1  (M)$ for some closed Riemannian manifold $M$, and $\{ M_i\}_{i\in\bbn}$ are the finite sheeted Galois covers of $M$ corresponding to $\{ N_i\}_{i\in\bbn}$, then we also have

\item There exists $\vare_5 > 0$ such that $h(M_i) \ge \vare_5$ (see Definition \ref{def1.7}) for every $i \in \bbn$.

    \item There exists $\vare_6 > 0$ such that $\lambda_1 (M_i) \ge \vare_6$ for every $i\in\bbn$ (see Proposition \ref{prop1.8}).

\end{enumerate}
\end{theorem}

The equivalence of all the properties except of $(v)$, can be found in [L1, Theorem 4.3.2].  For $(v)$ see [Sh1].

\chapter{Examples of Expanders}

It is by no means clear that expander graphs exist, though it is not difficult to prove their existence by random considerations.  Various deep mathematical tools have been used to give explicit constructions (Kazhdan property $T$, Ramanujan conjecture, etc.). We review these briefly here (again, more details are in [L1] and the references therein).  Most of the current chapter is devoted to a description of the two important developments of the last decade:

\begin{enumerate}

 \item All non-abelian finite simple groups are expanders in a uniform way, and

 \item Many linear groups have property $(\tau)$ with respect to congruence subgroups.
\end{enumerate}

The second result has important applications to number theory, group theory and geometry, to be described in later chapters.

\section{Random graphs and the Zig-Zag construction}

It is relatively easy to show that for a fixed $k\ge 3 $ there exists an $ \varepsilon  > 0$
s.t. a random $(n, k)$-graph, i.e. a $k$-regular graph on $n$ vertices, is an $\vare$-expander  with probability tending to 1 as $n$ goes to infinity.  The subtle issue (that we will ignore here) is how to describe a good model of random $(n, k)$-graphs. Anyway this has been known for many years (see \cite{L1} and \S 1.1 above).  More recently a much deeper result (Alon's conjecture) has been proved by Friedman [Fr].
\begin{theorem}
For every $\vare > 0$ and $k \ge 3$,
\be* Prob \big(\la(X) \le 2\sqrt{k-1} +\vare\big) = 1-o_k (1)\ee*
when $X$ is a random $(n, k)$-graph, i.e., almost every such $X$ is almost Ramanujan.
\end{theorem}

It is interesting to repeat here the open question mentioned in \S 1.2, whether for \emph{every} $k\ge 3 $ there exist infinitely many Ramanujan graphs.  While Theorem 2.1 hints toward a positive answer, it by no means implies that.  Moreover, one may be tempted to conjecture, following Theorem 2.1, that almost every $(n, k)$-graph is Ramanujan.  But quite a lot of computational data suggests to the contrary, namely that the probability of a $k$-regular graph on $n$ vertices to be Ramanujan tends, when $n\to \infty$, to  some constant strictly between 0 and 1. We do not know any result or even a conjecture that predicts what is this interesting number.

Of course, for any applications one wants explicit constructions.  A lot of work has been dedicated to this goal and various deep methods have been used.  In a breakthrough paper  Reingold, Vadhan and Wigderson \cite{RVW} showed that there is an elementary combinatorial way to build expanders via the Zig-Zag product of graphs, which they introduced.  We describe this subject here very briefly sending the reader to \cite{RVW} for details (or to \cite{HLW} for a very clear exposition).

The Zig-Zag product is a method that given two graphs $X$ and $Y$, where $X$ is an $(n, m)$ graph (i.e., $m$-regular graph on $n$ vertices) and $Y$ an $(m, d)$-graph, produces $X \circ Y$, an $(mn, d^2)$-graphs. (This is \emph{not} a commutative operation).    In [RVW], it is shown that one can bound the ``spectral gap" of $X \circ Y$ by the spectral gaps of $X$ and $Y$.  They then start for a small fixed integer $d$, with a $(d^4, d)$-graph $X=X_0$ with a good spectral gap (such a graph can be found by an exhaustive search in constant time) and define by induction $X_1 = X^2$ and $X_{n+1}  = (X_n)^2 \circ X$ for $n\ge 1$. (If $Y$ is a graph,  we denote $Y^2$ to be the graph with the same vertex set as $Y$, putting an edge between the end points of any path of length 2 in the original graph.) Note that $X_n$ is an $(d^{4n}, d^2)$-graph, so the family $\{ X_n\}^\infty_{n=1}$ is a family of $d^2$-regular graphs (independent on $n$). Induction and the spectral control give that this is a family of expanders.

This construction turns out to be quite useful in various applications in computer science, sometimes giving better results than using the expanders constructed by the other methods.  Still, as far as we know, it has not been used for applications in pure mathematics, which is our main interest in these notes.  For the kind of applications we emphasize here, one usually needs expanders which are Cayley graphs (or at least somewhat symmetric).  The Zig-Zag product has something to say about Cayley graphs of semi-direct products of groups (see \cite{ALW}, \cite{MW} and especially \cite{RSW}) but one still needs the other approaches.
These will be described in the next subsections.

\section{Kazhdan property $(T)$ and finite simple groups}

The seminal work of Kazhdan [Ka] on property $(T)$ of high-rank simple Lie groups and their lattices
(= discrete subgroups of finite covolume) opened the door for Margulis to give the first explicit examples of
expanders.  Let us repeat what we already said in \S 1.3.

\begin{prop}\label{ab} Let $\Gamma$ be a group with property $(T)$ generated by a finite symmetric set $\Sigma$, and let $\call = \{ N\big| N\triangleleft \Gamma, \; [\Gamma; N] < \infty\} $ be the family of finite index normal subgroups of $\Gamma$.  Then there exists an $\vare > 0$ such that all $Cay (\Gamma/N; \Sigma)$ are $\vare$-expanders.
\end{prop}

So, the issue (which is  non-trivial) is to show that such $\Gamma$'s exist.  This was done by Kazhdan and has been generalized substantially in recent years:

\begin{theorem}\label{abc}
Let $G$ be a simple Lie group (e.g. $G = SL_n(\bbr))$ of $\bbr$-rank $\ge 2$ (e.g. $n\ge 3$).  If $\Gamma$ is a lattice in $G$ (e.g. $\Gamma = SL_n(\bbz)$) then $\Gamma$ has property $(T)$.

\end{theorem}

Theorem \ref{abc} combined with Proposition \ref{ab} gives:
\begin{cor} For a fixed $n$, and a fixed finite symmetric set of generators $\Sigma$ of $SL_n(\bbz)$ the family of Cayley
graphs $Cay(SL_n (\bbz/p\bbz), \Sigma)$ are all $k$-regular $\vare$-expanders for $k = |\Sigma|$ and some
 $\vare = \vare(n, \Sigma) > 0$.
\end{cor}

$SL_n(\bbz/p\bbz)$ or more precisely $PSL_n(\bbz/p\bbz)$, the quotient by the center, is the prototype of finite
simple  groups of Lie type.  Many more infinite families of finite simple groups can be deduced, by a similar
method, to be expanders.  But a more challenging conjecture was put forward in 1989 by Babai-Kantor-Lubotzky [BKL]:
\begin{conj}\label{123} There exist $k\in\bbn$ and $\vare > 0$ such that every non-abelian finite simple  group $G$ has a symmetric
 set of generators $\Sigma$ of size $\le k$ such that $Cay (G;\Sigma)$ is an $\vare$-expander.
\end{conj}

For a number of years this conjecture has been open and even some suspicion arose (including by some of its proposers) that perhaps  it is not true and expansion is a property restricted to ``bounded rank". But it turned out that this conjecture is true and it is now fully proved and  its proof required an ensemble of very different methods.  The big breakthrough came with two works of Kassabov [K1] and  [K2] ([K1] was very much influenced by [KN] which in turn was modeled on [Sh3]).  Rather than being loyal to the historical development, let me  start with an even more recent result (which was influenced by [KN]).  To state the result we first need a definition.  Let $R$ be a ring with an identity.  For $n \ge 3$, define
$E_n(R)$ to be the multiplicative subgroup of the ring of $n\times n$ matrices $M_n(R)$ generated by $E_{ij} (r)$ for
all $1 \le i \neq j \le n$ and $r \in R$, where $E_{ij} (r)$ is the $n\times n$ matrix with $1$'s on the diagonal, $r$ at the $(i, j)$-entry and zero otherwise.
For many commutative rings, $E_n (R)$, at least for $n\ge 3$, is nothing more than $SL_n(R)$.  The groups $E_n (R)$ play an important role in algebraic $K$-theory.

If $R$ is a finitely generated ring then $E_n (R)$ is a finitely generated group (for $n\ge 3$).  We can now state:
\begin{theorem}[Ershov-Jaikin \cite{EJ}] Let $R$ be the  free ring $R=\break \bbz \langle x_1,\ldots, x_r\rangle$ in the non-commutative free variables $x_1,\ldots, x_r$.  Then for every $n\ge 3, E_n(R)$ has Kazhdan property $(T)$.
\end{theorem}

Kassabov's basic idea was to use such a result (well, he proved a weaker  version of it, which was sufficient) to deduce:

\begin{cor}\label{27} There exist $k\in\bbn$ and $\vare > 0$ such that for  every $n\in\bbn$ and every prime power $q\in\bbn$,
the group $SL_n(\bbf_q)$ has a set $\Sigma$ of $k$ generators for which
$Cay (SL_n (\bbf_q), \Sigma)$ is $\vare$-expander.
 Here, $\bbf_q$ is the field with $q$ elements.
\end{cor}

Indeed, take $R = \bbz \langle x_1, x_2\rangle$.  It is easy to see that $M_n(\bbf_q)$ is a (finite) quotient of $R$.  Thus by Theorem 2.6 and Proposition \ref{ab},
\be* Im \big(E_3 (R) \to E_3(M_n(\bbf_q)\big)\ee*
are expanders.  But this later group is
\be* E_3 \big(M_n (\bbf_q)\big)= SL_{3n} (\bbf_q).\ee*

    Now it is not difficult to deduce that $SL_n (\bbf_q)$ are expanders in a uniform way (i.e., with the same $k$ and  $\vare$ for all $q$ and for all $n\ge 3$) since $SL_{3n+1} (\bbf_q)$ and $SL_{3n+2} (\bbf_q)$ are bounded products of copies of $SL_{3n}(\bbf_q)$:
\begin{definition}\label{def2.1}  Let  $ \ga= \{ A_i \}_{i \in I}$ and $\calb = \{ B_j\}_{j\in J}$ be two families of finite groups.
We say that $\calb$ is a \emph{bounded product} of $\ga$ if there exists a constant $m\in\bbn$ such that for every
$B_j \in \calb$, there exist $A_{i_1,\ldots,} A_{i_m} $ in $\ga$ and homomorphisms $\varphi_{i_t} : A_{i_t} \to B$ such that
$B$ is the product (just as a set)  $\varphi_{i_1} (A_{i_1})\cdot\ldots\cdot \varphi_{i_m} (A_{i_m})$.
\end{definition}

The following easy Lemma is very useful:
\begin{lem}\label{29} In the notation above: if $\ga$ are expanders in a uniform way (i.e. same $k$ and $\vare$)  and $\calb$ is a bounded  product of $\ga$, then $\calb $ are also expanders in a uniform way (for some $k'$ and $\vare'$ which depend on $k$ and $\vare$).
\end{lem}

The next Theorem says that we can go much further than $SL_n$:
\begin{theorem}[Nikolov \cite{Ni}, Lubotzky \cite{L6}]\label{nikolov}
Let \be* \ga = \{ SL_n (\bbf_q) \big| n\ge 2,\;  q \ {\rm prime\ power}\} \ee* and $\calb$ is the family of all simple groups of Lie type excluding the Suzuki groups. The $\calb$ is a bounded product of $\ga$.
\end{theorem}

One can check in the proof that if one allows to use only $SL_n$ with $n\ge 3$, the result remains true for those groups in $\calb$ of rank $\ge 14$.  Thus Conjecture \ref{123} is valid for these groups. To handle all finite simple groups of Lie type one should handle $SL_2(\bbf_q)$. This will be done in 2.4 by a completely different method.  With the above results this will finish all groups of Lie type except the Suzuki.  They will be handled in 2.6 again by a completely different method.  But now we will handle first the case of Symmetric and Alternating groups which is of special interest.

\section{The symmetric groups}

Making the symmetric groups (or equivalently the Alternating groups) into a family of expanders in a uniform way has been a challenging problem for almost two decades, until it was solved by Kassabov [K2].

\begin{theorem}  There exists $k\in \bbn$ and $0 < \vare \in \bbr$ such that for every $n\ge 5, Sym(n)$ has a symmetric generating subset $\Sigma$ with $|\Sigma|\le k$ for which $Cay (Sym(n);\Sigma)$ is $\vare$-expander.
\end{theorem}
The same result holds also with $Alt (n)$-the alternating group instead of $Sym(n)$.  It suffices to prove the Theorem for one of these two cases.

Now, one can show that while $Alt(n)$ contains many copies of groups of Lie type it is not a bounded product of such groups, so the results of the previous section do not suffice.  Still Kassabov looked at $n$'s of the form $n=d^6$ for $d=2^{3r }-1$ for some $r\in\bbn$.  Based on ideas similar to the ones in the previous section, he shows that the groups $\Delta_r = SL_{3r} (\bbf_2)^{d^5}$ are $\vare_0$-expanders w.r.t. generating sets $F$ of size at most $40$.  He then embedded $\Delta_r$ in $Alt(n)$ in 6 different ways which give 6 copies of $F$ in $Alt(n)$. He then showed that $Alt(n)$ are uniformly expanders with respect to the union of these 6 sets.  It should be stressed that $Alt(n)$ is not a bounded product of these 6 copies and the argument is far more involved, working with the representation theoretic version of expansion.  One should work with the various irreducible representations of $Alt (n)$ and Kassabov divided them into two classes giving different arguments according to their Young diagrams. The reader is referred to [K2] for details and to [KLN] for a sketch of the proof.

\section{Property $(\tau)$, $SL_2$ and groups of low rank}

As already mentioned in \S 1.3, one does not need the full power of Property $(T)$ to deduce that the finite quotients of the finitely generated group $\Gamma$ give a family of expanders. Property $(\tau)$ (Definition  \ref{def1.6}) suffices.

The prototype of groups with $(\tau)$ w.r.t. some family is $\Gamma = SL_2(\bbz)$.

Let us set some notations:  $\Gamma = SL_2(\bbz)$ and for $m \in \bbn, \; \Gamma (m) = Ker (SL_2(\bbz) \to SL_2(\bbz/m\bbz))$ - the congruence subgroup $\mod m$.   The group $\Gamma$ acts on
$\bbh = \{ a + bi| a \in \bbr, 0 < b \in \bbr\}$ by Mobius transformations:
$\gamma = {a b \choose c d} \in\Gamma$ and $z \in\bbh$, then $\gamma (z) = \frac{az+b}{cz+d}$. The upper half plane  $\bbh$ is endowed with a Riemannian metric of constant curvature $-1$.  This is the Hyperbolic plane.  The quotients $\Gamma(m) \setminus \bbh$ are (non-compact) Riemann surfaces of finite volume.

\begin{theorem}[Selberg\cite{Sel}]\label{se} For every $m\in \bbn$,  $\lambda_1 \big(\Gamma(m)\setminus\bbh\big) \ge \frac{3}{16}$.
\end{theorem}
Selberg conjectured that $\frac 14$ is the right lower bound.  The current world record is $\lambda_1\ge 0.238$  due to Kim and Sarnak \cite{Ki}.

Anyway Theorem \ref{14}
  gives:
 \begin{cor}  For a fixed set of generators $\Sigma$ of $\Gamma = SL_2 (\bbz)$, e.g., $\Sigma = \big\{ {1 \; \pm 1 \choose 0  \; \; \; 1}, {1 \; \; \; 0\choose \pm 1 \; \;1} \big\}$, the Cayley graphs $Cay\big( SL_2(\bbz/m\bbz);\Sigma\big)$ are all $\vare$-expanders for some $\vare > 0$ which may depend on $\Sigma$ but not on $m$.
  \end{cor}

It should be stressed that $\Gamma = SL_2(\bbz)$ has negative solutions to the congruence subgroup problem and it has many more finite quotients than just $SL_2(\bbz/m\bbz)$.  The family of \emph{all} finite quotients of $\Gamma$ do not form a family of expanders.  So in the terminology of Definition \ref{def1.6}  above, $\Gamma$ does not have property $(\tau)$, but it has property $(\tau)$ w.r.t. the family of congruence subgroups.

The last Corollary implies in particular that the groups\break $\{ SL_2(p)|p \; \mathrm{prime}\}$ can be made into a family of 4-regular Cayley graphs which are expanders uniformly (i.e. same $\vare$).

Analogous results for arithmetic groups in positive characteristic such as $SL_2(\bbf_p [t])$ or $ SL_2(\bbf_p [t, t^{-1}])$ (when this time a result of Drinfeld replaces the Theorem of Selberg) can make, \emph{for a fixed $p$}, the family $\{ SL_2 (\bbf_{p^\a}) \big| \a \in\bbn\}$ into a family of expanders.

Can we make all of these families together
$\{ SL_2 (\bbf_{p^{\a}}) \big| p \; { \rm prime}, \a \in \bbn\}$ into a family of expanders?

The answer is yes - but the proof is more subtle.  See [L6] (and a sketch in [KLN]).  Here we only mention that the proof uses the \emph{explicit} constructions of Ramanujan graphs (as a special case of Ramanujan complexes) in [LSV2]  (see also [LSV1]).  It is shown there that for a fixed $p$, $SL_2(\bbf_{p^{\a}})$ are  $(p+1)$-regular Ramanujan graphs.  The $p+1$ generators involved are the conjugates of a fixed element $C_{p,\a}$ of $SL_2(\bbf_{p^\a})$ by a fixed non-split torus $T\subseteq SL_2(\bbf_p) \subseteq SL_2 (\bbf_{p^\a})$.
Then use the fact that $SL_2(p)$ are expanders with respect to $\Sigma_0 = \big\{ {\; 1 \; \; \, \pm 1\choose 0 \; \; \;\;  1} ,
{\; 1 \; \; \; \; \; 1\choose \pm 1 \; \; \;\;  1}\big\}$ to deduce that $SL_2 (\bbf_{p^\a})$ are expanders w.r.t.  the symmetric set of six generators $\Sigma_0 \bigcup \{ C^{\pm 1}_{p, \a}\}$ in a uniform way.  So  Selberg Theorem and Drinfeld solution to the positive characteristic Ramanujan conjecture for $GL_2$ are both needed for that goal.

Once all the $SL_2(\bbf_{p^\a})$ are expanders and so all $\{ SL_n (\bbf_{p^\a})\big | p, n, \a \in \bbn, \; p \; \rm{prime}\}$ are expanders
in a uniform way (same $k$, same $\vare$).  We can appeal again to Theorem \ref{nikolov}  and Lemma \ref{29} to deduce that all finite simple groups of Lie type, with the possible exceptions of the Suzuki groups, are expanders in a uniform way.

Suzuki groups have to be excluded here as they do not contain copies of $SL_n(\bbf_{p^\a})$. Indeed their order is not divisible by 3.  (A classical result of Glauberman  at the early days of the classification project of finite simple groups asserts that this property characterizes them! See [Gl]). But the Suzuki groups are not exceptional for our problem - i.e. they are also expanders.  But this requires another method and has to wait for \S 2.7.

\section{Property $(\tau)$ with respect to congruence subgroups}

Selberg Theorem \ref{se} was a starting point for many works which extended it to general arithmetic groups.  The results are of importance in number theory (automorphic forms), representation theory and geometry.  In this section, we will describe them from our perspective.

Let $k$ be a global field, i.e., a finite extension of $\bbq$ or of $\bbf_p(t)$.  Let $G$ be a simple algebraic group defined over $k$ with a fixed embedding $\rho: G \hookrightarrow GL_n$ for some $n$.  Let $\theta$ be the ring of integers of $k$ and $S$ a finite set of valuations of $k$ containing $S_\infty$ - the set of archimedean valuations.  Let $\theta_S = \{ x \in k\big| v(x) \ge 0,  \;  \forall v \notin S\}$ - the ring of $S$-integers, so $\theta_S = \theta$ if $S = S_\infty$.  Let $\Gamma = \rho \big(G(k)\big) \cap GL_n(\theta_S)$.  A subgroup of $G$ commensurable with $\Gamma $ is called an $S$-{\it arithmetic subgroup} of $G$.  For a non-zero ideal $I$ of $\theta_S$ (which is always of finite index) we denote $\Gamma(I) =\break Ker \big(\Gamma \to GL_n(\theta_S/I)\big)$.
An $S$-arithmetic subgroup of $G$ containing $\Gamma (I)$ for some $I$ is called a $(S-)$ {\it congruence subgroup}.

While the definition of $\Gamma$ may depend on the choice of the representation $\rho$ the classes of arithmetic and congruence subgroups do not.
\begin{definition}\label{214}
We say that $\Gamma$ has the Selberg property if it has property $(\tau)$ with respect to the congruence subgroups $\{ \Gamma (I)\}_{0 \neq I \tle \theta_S}$
\end{definition}

Again, if true for $\Gamma$, then it is true for all the arithmetic groups in its commensurability class.

The group $\Gamma = G(\theta_S)$ sits as an irreducible lattice in the Lie group $H=\mathop{\Pi}\limits_{v\in S} G(k_v)$ where $k_v$ is the completion of $k$ w.r.t. $v$.

Recall that a \emph{lattice} $\Lambda$ in $H$ is a discrete subgroup where $\Lambda \setminus H$ carries an $H$-invariant finite measure.  It is \emph{irreducible} if its projection to each $G(k_v)$ is dense. In many cases $H$ has Kazhdan property $(T)$:

\begin{theorem}[Kazhdan]\label{kazhdan} If $k_v$-rank $(G) \ge 2$ for every $v\in S$, then $H = \underset{v\in S}{\Pi} G(k_v)$ and all its lattices have property $(T)$.
\end{theorem}

For having $(\tau)$ we need less:
\begin{theorem}[Lubotzky-Zimmer \cite{LZi}] If one of the non-compact factors of $H$ has $(T)$, then all irreducible lattices have property $(\tau)$.
\end{theorem}

Selberg Theorem \ref{se}  shows that $\Gamma = SL_2 (\bbz)$ which has neither $(T)$ nor $(\tau)$, still has $(\tau)$ w.r.t. congruence subgroups. This has been extended to all arithmetic groups.  This is the work of many people.  The most general method  (and actually also the simplest!) is due to Burger and Sarnak ([BS]) who proved:

\begin{theorem}\label{217} If $L_1 \le L_2$ are two non-compact simple Lie groups with arithmetic lattices $\Lambda_i \le L_i, i = 1, 2$ and $\Lambda_1 = L_1 \cap \Lambda_2$.  Then:
\begin{enumerate}[(i)]
\item If $\Lambda_1$ has property $(\tau)$ so does $\Lambda_2$.
    \item If $\Lambda_1$ has the Selberg property, so does $\Lambda_2$.
\end{enumerate}
\end{theorem}

Many (in some sense ``most") simple $k$-algebraic groups $G$ contain a copy of $SL_2$, so Theorems \ref{se} and \ref{217}  imply the Selberg property for the arithmetic subgroups of $G$. By Galois cohomology method one can classify the arithmetic lattices for which this method does not apply.  These need some other (more difficult) techniques.  This was done by [Cl] using automorphic forms methods.  As a result it is now known that all arithmetic lattices in semi-simple groups over local fields of characteristic zero have the Selberg property. As far as we know this has not been completed yet for the positive characteristic case.

\section{Sum-products in finite fields and expanders}

The results described in the previous section gave a fairly complete picture on the congruence quotients of an arithmetic group of the form $\Gamma = G(\theta_S)$ described there as being expanders with respect to generators coming from ``the mother group" $\Gamma$.
For example, the family
\be*
\{ Cay\big(SL_2(\bbf_p);
{\textstyle{{1 \; \; \; \pm 1\choose 0 \; \; \; 1}, {1\; \; \; 0\choose \pm 1 \; \; \; 1}}}
\big) \big| p \ {\rm prime}\}\ee*
form a family of $\vare$-expanders for some $\vare > 0$ since
${1 \; \; \; 1\choose 0 \; \; \; 1}$
and $  {1\; \; \; 0\choose  1 \; \; \; 1}$ generate $SL_2(\bbz)$.  A similar conclusion (for a different $\vare > 0$) is true for the family
\be*  \{ Cay\big(SL_2(\bbf_p);
{\textstyle{{1 \; \; \; \pm 2\choose 0 \; \; \; 1}, {1\; \; \; 0\choose \pm 2 \; \; \; 1}}}
\big) \big| p> 2 \ {\rm prime}\}\ee*
even though $\Gamma$ is not generated by
${1 \; \; \; 2\choose 0 \; \; \; 1}$
 and
${1\; \; \; 0\choose 2 \; \; \; 1};$
 these two matrices generate a finite index subgroup of $\Gamma$ and essentially the same arguments as before applied also for it.

But now what about
\be*\{ Cay \big(SL_2(\bbf_p);
{\textstyle{{1 \; \; \; \pm 3\choose 0 \; \; \; 1}, {1 \; \; \; 0\choose \pm 3 \; \; \; 1}}}
\big) \big | \; p > 3 \;  {\rm prime} \} ? \ee*

Is this a family of $\vare$-expanders for some $\vare > 0$?  The issue here is that the subgroup $\Lambda = \langle {1 \; \; \; 3\choose 0 \; \; \; 1}, {1 \; \; \; 0\choose 3 \; \; \; 1} \rangle$ is of infinite index in $SL_2(\bbz)$, still when taken $\mod p$, the image of $\Lambda$ generates $SL_2(\bbf_p)$ for every $p > 3$.  Combinatorially one should expect a similar $\vare$-expander conclusion for them but the methods of the previous sections do not apply here. This problem was presented in 1992 in [L2] and was popularized under the nickname (given by Alex Gamburd) ``Lubotzky 1-2-3 problem".

In fact this 1-2-3 problem is just an attractive special case of a much more general problem: Let $\Gamma = G(\theta_S)$ as in the previous section, but for simplicity of notation, let's take $k = \bbq, \; \theta  = \bbz$ and $S = S_\infty$, so $\Gamma = G(\bbz)$, e.g. $\Gamma = SL_d(\bbz)$.  Let $\Lambda$ be a finitely  generated subgroup of $\Gamma$, generated by a set $\Sigma$, which is Zariski dense in $\Gamma$.  Note that being dense in the Zariski topology is quite a weak assumption.  For example for $\Gamma = SL_2(\bbz)$, every subgroup $\Lambda$ which is not virtually cyclic, e.g. $\Lambda = \langle {1 \; \; \; 3\choose 0 \; \; \; 1}, {1 \; \; \; 0\choose 3 \; \; \; 1}\rangle$ is Zariski dense.  Assume further that $G$, as an algebraic group is connected, simply connected and simple (e.g. $G = SL_d$).  Then the  \emph{strong approximation theorem for linear groups} (\cite{MVW}, \cite{No}, \cite{W}, \cite{Pi} - see [LS, Window 9] for an exposition) says that there exists $m_0 \in\bbn$ such that for every $m\in\bbn$ with $(m, m_0) = 1$ the projection of $\Lambda$ to $G(\bbz/m\bbz)$ is onto.  In other words it says that $Cay \big( G(\bbz/m\bbz); \Sigma\big)$ is a connected graph.  Being expander is a very strong form of being connected.  It thus naturally suggests the question whether these graphs form a family of $\vare$-expanders.
The so-called ``Lubotzky 1-2-3 problem" is just a baby-version of this much more general question.

First steps and several interesting partial results toward the 1-2-3 problem were taken in \cite{Ga1} and \cite{Sh2}.  But the main breakthroughs came in recent years, starting with the work of Helfgott and continued by others.  We  now turn to describe these developments.

Let us start by stating the main result of [H]:
\begin{theorem}\label{218}  Let $G = SL_2(\bbf_p)$ and $A$ a generating subset of $G$.  Let $0 < \delta < 1$ be a constant. Then
\begin{enumerate}[(a)]
\item If $|A| < |G|^{1-\delta}$, then $|A \cdot A \cdot A|\ge C |A|^{1+\vare}$ where $C$ and $\vare$ depend only on $\delta$.

    \item If $|A| \ge |G|^{1-\delta}$ then $ A\cdot\ldots\cdot A=G$, i.e.,  the product of $k$  copies of $A$ is $G$, where $k$ depends only on $\delta$.

        \end{enumerate}
        \end{theorem}

        Before elaborating on its importance for expanders, let us put it in a more general context.

        The sum-product results form a body of various theorems asserting that if $F =\bbf_p$ is a finite field of a prime order $p$ and $A$ a subset of $\bbf_p$, which is not too large then either the set of products $A\cdot A = \{ a\cdot b|a, b\in A\}$ or the set of sums $A+A = \{ a+b| a, b \in A\}$ is significantly larger than $A$.  Here is a typical result in this area, called also ``additive combinatorics" (\cite{TV}).

\begin{theorem}[{[BKT]}]\label{219} If $A$ is a subset of $\bbf_p$, $p$ prime, with $p^\delta \le |A| \le p^{1-\delta}$ for some $\delta > 0$, then $|A + A| + |A\cdot A|\ge c|A|^{1+\vare}$, where $c$ and $\vare$ depend only on $\delta$.

\end{theorem}

The main idea of Helfgott was to convert the growth of a subset $B$ of $SL_2(\bbf_p)$ when taking product $B\cdot B\cdot B$ with the growth of $A = tr (B) = \{ tr(g) \big| g \in B\}$ as a subset of $\bbf_p$ under sums and products.  He also showed that the sizes of $B$ and $A$ can teach a lot about each other. This enabled him to deduce Theorem \ref{218} from Theorem \ref{219}. His work is quite complicated from a technical point of view.  Some subsequent works simplified and extended his work - see below - and eventually made the conclusion free of the use of sum-products results.  Still, various ideas of Helfgott  are still crucial also in those extensions.

An interesting corollary of Theorem \ref{218} is that there exists a constant $C$ such that for every set of generators $\Sigma$ of $SL_2(\bbf_p)$,
\be* diam (Cay(SL_2(\bbf_p); \Sigma)) \le \log (p)^C.\ee*
This was the first infinite class of groups for which the following long-standing conjecture of Babai was proved:
\begin{conj}\label{220}
There exists a constant $C$, possibly $C = 2$, such that for every finite simple group $G$ and for every set of generators $\Sigma$, the diameter is polylogarithmic  (i.e., diam $Cay(G;\Sigma)) = O((\log|G|)^C)$.
\end{conj}

The example $Cay(Sym(n); \tau = (1,2), \sigma^\pm= (1,2,\ldots, n)^{\pm 1}\}$ and similar ones for $Alt (n)$ show that one cannot expect better than $C=2$ (see \cite{L1}).

While Helfgott's result solved Babai's conjecture for $SL_2(\bbf_p)$, it fell short of showing that these are expanders.  (By the way, expanders give rise to logarithmic diameter -  i.e. $C=1$ in the last conjecture). It did not solve the 1-2-3 problem either. But shortly afterwards Bourgain and Gamburd \cite{BG1} made a second major breakthrough establishing the desired expansion by introducing their fundamental flattening lemma technique and coupling it with more standard techniques from the representation theory of these finite simple groups.
\begin{theorem}\label{221}  For any $0< \delta \in\bbr$ there exists $ \vare = \vare (\delta) \in\bbr$ such that for every prime $p$, if $\Sigma$ is a set of generators of $SL_2(\bbf_p)$ such that girth $ \big(Cay (SL_2(\bbf_p); \Sigma)\big) \ge \delta\log p$ then $Cay \big( SL_2(\bbf_p); \Sigma\big)$ is an $\vare$-expander.
\end{theorem}

 The girth of a graph is the length of the shortest non-trivial closed path in the graph.

This theorem solves, in particular, the 1-2-3 problem: an  easy argument (going back to [M1]) shows that \be* {\rm girth\ } \big(Cay(SL_2(\bbf_p);
{\textstyle{{\; \;1 \; \; \; \pm 3\choose 0 \; \;  1}, {1\; \; \;\;\, 0\choose \pm 3 \; \; \; 1}}}
)\big)\ge \delta\log p\ee*
for some $0 < \delta $ independent of $p$.  Moreover, if $\Sigma$ is  a free set of generators of a non-abelian free subgroup of $SL_2(\bbz)$, then girth $\big(Cay(SL_2(\bbf_p); \Sigma)\big)$
is logarithmic in $p$ and hence these are expanders.  In fact, the last conclusion holds for every Zariski dense subgroup $\Lambda$ of $SL_2(\bbz)$ as every such subgroup contains a non-abelian free group.  This is easy for $SL_2(\bbz)$ but true also for the more general case of $\Lambda$, which concerns us, by a well known result of Tits [Ti].  It actually implies that we can assume $\Lambda$ is a non-abelian free group.

The reader may note that we have stopped discussing  general congruence quotients $SL_2 (\bbz/m \bbz)$ for $m \in \bbn$ and stuck to $m = p$ a prime $SL_2(\bbz/p\bbz) = SL_2(\bbf_p)$.
Well, the extension to general $m$ required more effort. It was first done in \cite{BGS2} for natural numbers $m$ which are square free and eventually for all $m \in\bbn$ in [BV].  We will come back to this issue later.  The case of $m$ square free is especially important for the sieve methods in Chapters 4 and 5.

The dramatic breakthroughs have continued even further; first in parallel by two groups of researchers Breuillard-Green-Tao ([BGT1], [BGT2]) and Pyber-Sza\'bo ([PS1], [PS2]) and secondly by Varju \cite{V} (followed by \cite{SGV}).
The first two groups proved essentially the same result (there are some differences but for our impressionistic picture we can ignore them).
\begin{theorem}\label{222} Let $r \in\bbn$ be fixed.  Then for every finite simple group $G$ of Lie type of rank at most $r$ and for every subset $A$ of $G$ which generates  $G$, either $A\cdot A\cdot A = G$ or $|A\cdot A\cdot A|\ge |A|^{1+\vare}$ where  $\vare$ depends only on $r$.\end{theorem}

The reader may check that this generalizes Helfgott's result (Theorem \ref{218}) from $SL_2(\bbf_p)$ to all finite simple groups of Lie type \emph{of bounded rank}.  It should be mentioned that shortly after Theorem \ref{218} was proved, it was shown by Nikolov and Pyber [NiPy] that part (b) of that theorem follows quickly from a result of Gowers [Go] and in more general situation, i.e., finite simple groups of Lie rank at most $r$ with $k=3$ for some $\delta = \delta(r)$ depending only on $r$.  (See \cite{BNP} for more).  This handles the case of ``large subsets" and  the main novelty of Theorem \ref{222} is the case when $|A| < |G|^{1-\delta(r)}$.

Theorem \ref{222} extends Helfgott's result from $SL_2$ to any bounded rank finite simple group. In particular it proves the Babai Conjecture for this case, namely:
\begin{cor}\label{223}
For $r\in\bbn$, there is a constant $C=C(r)$ such that for every finite simple group of Lie type $G$ of rank at most $r$ and every symmetric set of generators $\Sigma$ of $G$,
\be* diam \big(Cay(G;\Sigma)\big) \le (\log|G|)^C.\ee*
\end{cor}

The conjecture is still open for the unbounded rank case.
It should be stressed  that in Theorem 2.22, $\varepsilon$ does depend on $r$.  In fact, it is shown in \cite{PS2} that $\vare = O(\frac 1r)$.  Still one can expect that $C$ of Corollary 2.23 is independent of $r$.
For the bounded rank case one may conjecture that the right bound is $C(r)(\log |G|)$ rather than $(\log|G|)^{C(r)}$ - see more in \S 2.7.

Helfgott's result for $SL_2$ led to  Bourgain-Gamburd Theorem \ref{221}.  It naturally suggests to expect a similar result for bounded rank simple groups.  The following theorem is the second breakthrough (proved first in \cite{V} for $SL_n$ and later in \cite{SGV} in general). It is in one way weaker and in another way stronger than the expected analogue.

\begin{theorem}[A. Salehi Golsefidy-P. Varju \cite{SGV}]\label{224} Let $\Gamma$ be a finitely generated subgroup of $GL_n(\bbq)$, so $\Gamma \subseteq GL_n(\bbz_S)$ for some finite set of primes $S$.  Let $m_0 = \underset {p \in S}{\Pi p}$.  For every $q$ prime to $m_0$, let $\Gamma (q) = Ker\big( \Gamma\to GL_n(\bbz_S/q\bbz_S))$.
Then $\Gamma$ has property $(\tau)$ with respect  to the family $\{ \Gamma (q) \big| q\; {\rm is\ square\ free\ }\}$ iff $H^0$, the connected component of the Zariski closure $H = \ol{\Gamma}^z$, is perfect (i.e., does not have an abelian quotient).
\end{theorem}

The only if  part is easy. The main point is the if part. It does not exactly generalize Theorem 2.21 of Bourgain-Gamburd, but rather generalizes Bourgain-Gamburd-Sarnak \cite{BGS2} which gives expansion for  to square-free congruence quotients but only for $SL_2$.   As we will see in the next chapters, this is an extremely useful result with important applications to number theory, group theory and even geometry.  The reader may note that even though we formulated the result  over $\bbq$, it holds over any number field by restriction of scalars and for most applications one can reduce anyway to this case.

It will be desirable to know the above result for $\Gamma (q)$ for all $q$'s without the restriction to squarefree (see also \cite{BG3} and \cite{BG4}),  though at this point we do not see applications to this more general statement.  So far this was proved only for $\Gamma = SL_n(\bbz)$ by Bourgain and Varju (\cite{BV} using \cite{BFLM}).  One can fantasize on a much more general statement, which will be the ultimate generalization of the 1-2-3 problem:
\begin{conj}\label{225}  Let $\Gamma$ be a finitely generated subgroup of $GL_n(F), F$ a field.  So, $\Gamma\subset GL_n(R), \; R$ a finitely generated domain. Assume $H^0$ is perfect where $H$ is the Zariski closure of $\Gamma$ and $H^0$ its connected component.  Then $\Gamma$ has $(\tau)$ w.r.t. to the family $\Gamma (I) = Ker \big(\Gamma \to GL_n(R/I)\big)$ when $I$ runs over all the finite index ideals of $R$.
\end{conj}

A proof of the positive characteristic part of this conjecture will have new applications.

In light of Theorem 2.6 one can even speculate on  a more general version of conjecture \ref{225} but this conjecture is general enough to cover any applications in sight.

\section{Random generators and worst case generators}

In Sections 2.2, 2.3 and 2.4 we described how all the non-abelian finite simple groups, except for the Suzuki groups, can be made into a family of expanders uniformly.  We showed that there exists $k \in\bbn$ and $0 < \vare \in\bbr$ such that for every such $G$ there exists an explicitly given symmetric subset $\Sigma$ of $G$ of size at most $k$ such that $Cay (G; \Sigma)$ is an $\vare$-expander. We did not bother to write the sets $\Sigma$ explicitly but the method was explicit and if one wants, such a $\Sigma$  can be presented.  Now, for the Suzuki groups such $\Sigma$ exists but in a non-explicit way.

\begin{theorem}[Breuillard-Green-Tao {[BGT3]}]\label{th2.26} There exists an $0 < \vare \in\bbr$ such that for every Suzuki group $G=Sz (2^{2\ell + 1})$, for almost every pair of elements $(x, y) \in G\times G$, the Cayley graph $Cay(G; \{ x^{\pm 1}, y^{\pm 1}\})$ is an $\vare$-expander.
\end{theorem}

So altogether Conjecture 2.5 is now a Theorem!  This result handles the ``\emph{best-case scenario}", i.e., there exists a set of generators $\Sigma$ of $G$ with the desired property.  What about random set?

Recall the well known result:
\begin{theorem}[\cite{D}, \cite{KL}, \cite{LiSh}] Two random elements of a finite simple group $G$ generate $G$ with probability going to 1 when $|G|$ is going to infinity.
\end{theorem}

Another way to state the last result is that a random pair of elements gives rise to a connected Cayley graph.  Is this graph expander?  In a more precise formation.

\begin{open} \label{2.28}  Is there an $0 < \vare \in\bbr$, such that $Prob (Cay (G;\{ x^{\pm 1}, y^{\pm 1} \})$ is $\vare$-expander) is going to one when $G$ runs over the non-abelian finite simple groups with $|G| \to \infty$ and $x$ and $y$ are chosen randomly and uniformly from $G$?
\end{open}

It was recently proved by Breuillard, Green, Guralnick and Tao (\cite{BGGT2}) that the answer to this problem is yes, if one restricts himself to groups of bounded Lie rank.  This of course generalizes Theorem 2.26 as well as a similar result which was known before for $SL_2$ (\cite{BG1} and \cite{Di} using \cite{GHSSV}). It can replace \cite{L6} as a proof for the bounded case of Conjecture 2.5. (The proof in \cite{L6} used deep results from automorphic forms like Selberg theorem and Drinfeld solution to the characteristic $p$ Ramanujan conjecture - but gave explicit generators).

One can ask for even more: Is it possible that there exists $\vare > 0$ such that $Cay(G;\{x^{\pm 1}, y^{\pm 1} \})$ is $\vare$-expander for {\it every} choice of generating set $\{ x, y\}$ for $G$ and any  non-abelian finite simple group (``worse case scenario"; compare to Babai conjecture 2.20).  In the general case the answer is certainly - no! The family $Alt(n)$ has generators which do not give rise to a family of expanders. (For $Sym(n)$ one can take $\{\tau=(1,2), \sigma = (1, 2, \ldots, n)\}$:  $Cay(Sym(n); \{ \tau, \sigma^{\pm 1}\} )$ are not expanders - see \cite[Example 4.3.3(c)]{L1} - and from this one can deduce a similar result for $Alt(n)$). It seems likely that for a family of finite simple groups of unbounded rank, one can always find    ``worse case generators" which are not expanders.  But one may suggest:

\begin{conj}\label{co2.29} For every $r \in \bbn$ there exists $\vare = \vare(r)$ such that $Cay(G;\{x^{\pm 1}, y^{\pm 1} \})$ is an $\vare$-expander for {\it every}  finite simple group $G$ of Lie type and rank at most $r$ and for {\it every} generating set $\{ x, y\}$ of $G$.
\end{conj}

One may want to compare Conjecture 2.29 with Corollary \ref{223} which gives  a weaker statement.  An intermediate step would be to prove that
\be* diam \; Cay(G;  \{x^{\pm 1}, y^{\pm 1} \}) = O_r(\log |G|)\ee* where the implied constant depends only on $r$.

As of now the only result concerning Conjecture \ref{co2.29} is:
\begin{theorem}[Breuillard-Gamburd \cite{BGa}]\label{th2.30}  There exists $0 < \vare \in\bbr$ and an infinite set of primes $\calp$ such that for every $p\in\calp$ and \emph{every } generating set $\{ x, y\}$ of $SL_2(\bbf_p), \; Cay(SL_2(\bbf_p); \{x^{\pm 1}, y^{\pm 1} \})$ is an $\vare$-expander.
\end{theorem}

Their interesting method is not explicit.  They prove the existence of such a set $\calp$ by a non effective method.

\chapter{Applications to computer science}

Expander graphs play an important role in computer science with numerous
applications in many subareas. They appear as basic building blocks
at various networks, give error correcting codes, are used for derandomization
of various probabilistic algorithms and more. Many of the applications
are presented in \cite{HLW} and the reader is encouraged to consult
them, either in details or at least to get impressed by the wide spectrum
of applications. We chose to present one ``real" application (to error correcting codes) and one theoretical application to the theory of computation (the analysis of the product replacement algorithm).

\section{Error correcting codes}

Error collecting codes is a collective name for various methods that enable sending  messages of information through  noisy channels. The most common model deals with  sending a block of $k$ bits of information, i.e. a vector $v$ in ${\bold F}^k_2$.  Instead of $v$, one sends $Tv\in {\bold F}^n_2$ when $n > k$, i.e. a longer vector, but with the hope that if the noise will cause $t$ mistakes (i.e. switching $0$ to $1$ or vice versa, in $t$ coordinates) the receiver will be able to correct it back to the right vector. This can happen if for every $v_1\neq v_2\in {\bold F}^k_2$, $\text{dist}(Tv_1,Tv_2) > 2t$ where for $x,y\in {\bold F}^n_2$, we take $\text{dist}(x,y) =$ the number of bits in which they are different.

It is usually convenient to use a linear transformation for $T$, in which case \newline
 $C:= T({\bold F}^k_2)$ is a linear subspace of ${\bold F}^n_2$.  Moreover, as $\text{dist}(x,y) = \text{dist} (x-y,\vec 0)$, this code can correct $\lfloor\frac{d-1}{2}\rfloor$ errors, when $d= d(C) = \text{min}\{ \text{dist} (x, \vec 0)\mid \vec 0\neq x\in C\}$.  This leads to define:

\begin{definition}  A $(n,k,d)$-code is a linear subspace $C$ of ${\bold F}^n_2$ of dimension $k$ of distance $d(C) = d$.
\end{definition}

In coding theory,  one is interested in  ``good codes'', i.e., a family of $(n,k,d)$-codes with dimension $k$ and distance $d$ both growing linearly with $n$.  So let us denote
$r(C) =\frac{k}{n}$ and $\delta (C) = \frac{d}{n}$, the rate and the relative distance of the code $C$. So, a family of codes, of dimensions going to infinity, is good if there exists $\varepsilon > 0$ such that $r(C)$ and $\delta (C)$ are both at least $\varepsilon$.

The subspace $C$ is defined by linear equations. The code (or more precisely the family of codes) is called LDPC (low density parity check) if for some fixed constant $\ell$, all these linear equations are $\ell$-sparse, i.e. each such equation touches at most $\ell$ variables. Another way to say this is that the ``parity check'' matrix $H$ defining $C$, i.e. the $(n-k)\times n$ matrix $H$ with $C = \{ x\in {\bold F}^n_2 \mid Hx = \vec 0\}$, has at most $\ell$ non-zero entries in each of its rows. (Of course, such $H$ is not unique; we say that $C$ is LDPC if such an $H$ exists).

It has been known for a long time that LDPC good codes do exist.  This was first shown by random considerations (see \cite{HLW} and the reference therein). In 1996,  Sipser and Spielman \cite{SS} gave an explicit construction based on expander graphs.

To describe their work, let us start  with a simpler construction of codes based on graphs (sometimes called ``cycle graph codes'') as follows:

Let $X = (V,E)$ be a connected $r$-regular graph on $m$ vertices. So $\vert E\vert = \frac{mr}{2}$. Let ${\bold F} = {\bold F}_2$ and ${\bold F}^E$ be the space of functions from $E$ to ${\bold F}$. This can be thought of as the ${\bold F}$-vector space with basis $E$ or also  as the set of all subsets of $E$ (where $A+B = (A\cup B)\backslash (A\cap B)$). Let $C$ be the subspace of ${\bold F}^E$ spanned by all cycles of $X$. A simple argument shows that if we consider ${\bold F}^E$ as the
space of functions, then $f\in C$ iff for every $v\in V$,
\begin{equation}\tag{${3.1}$}
\sum\limits_{v\in e\in E} f(e) = 0.
\end{equation}

(In fact, this is the same argument which enabled Euler to prove that there is no Eulerian path in Konegsbourg; i.e. there is such a path iff the degree of every vertex is even). This last remark shows that $\dim (C)\leq \vert E\vert - \vert V\vert$.  But, in fact, the sum of all the defining equations is $0$ and one can prove that $\dim(C) = \vert E\vert - \vert V\vert + 1$. Note that each one of the defining equations (3.1) has support exactly $r$, so if we take a family of $r$-regular graphs we get a family of LDPS codes with  rate $= 1 - \frac{2}{r}$. But, unfortunately, for these codes, the distance is logarithmic in the dimension rather than linear. Indeed, it is easy to see that the distance of this code is exactly girth$(X)$, the girth of the graph $X$, i.e., the length of the shortest non-trivial closed cycle in $X$. By the well known and easy Moore inequality, for an $r$-regular graph $X$ on $n$ vertices, $\text{girth} (X)\leq 2 \log_{r-1}(n)$, which implies that the cycle code cannot be good.

To overcome this, Sipser and Spielman used the following idea (which in some sense goes back to Tanner [Ta]); Choose a ``small'' code $C_0$ inside ${\bold F}^r$ with  rate $r_0$ and relative distance $\delta_0$. For every $v\in V$ give the edges coming out of $v$, labels $1,\ldots, r$, and denote them as $e_v(1),\ldots ,e_v(k)$ (we do not require any compatibility here: the same edge can be labeled $i$ when it comes out of $v$ and $j$ when it comes out of $w$). We now define $C(X,C_0)$ to be the subspace of all functions $f\in{\bold F}^E$ such that for every $v\in V$, $(f(e_v(1)),\ldots , f(e_v(k)))$ is a vector in $C_0$. Namely, these are the functions which are ``locally'' in $C_0$, i.e. what every vertex ``sees'' in its star is a vector of $C_0$.

\begin{theorem}[Sipser-Spielman {[SiSp]}]    The code $C(X,C_0)$ has relative rate at least $2r_0-1$ and relative distance at least $(\frac{\delta_0-\lambda}{1-\lambda})^2$ where $\lambda = \lambda (X) = \frac 1r\text{max}\{ \mu \mid \mu ~ \text{an eigenvalue of $X$}, \mu\neq r\}$.
\end{theorem}

The theorem gives an explicit construction of LDPC good codes. Indeed, let $X$ be an  $r$-regular Ramanujan graph, so $\lambda (X) \leq \frac{2\sqrt{r-1}}{r} \leq \frac{2}{\sqrt{r}}$. Pick a code $C_0$ in ${\bold F}^r$ with rate $> \frac{1}{2}$ and relative distance $> \frac{2}{\sqrt{r}}$. Such codes do exist as can be seen by either random consideration (and as $r$ is fixed, we are allowed to pick one randomly) or by one of the many classical methods (note that we only ask the relative distance to be more than $\frac{2}{\sqrt{r}}$, so it does not have to be ``good''). Theorem 3.1 now ensures that $C(X, C_0)$ is good.  Finally, the code
$C(X,C_0)$ is LDPC since every defining equation touches only the $r$ edges adjacent to a vertex $v$ (the equations which force it to be in $C_0)$.

The proof of Theorem 3.1 is not difficult. As $C_0$ is defined by $(1-r_0)r$ equations, $C=C(X,C_0)$ is defined by $(1-r_0)rm$ equations on $\vert F\vert = \frac{rm}{2}$ variables, so $\dim C\geq \frac{rm}{2} - (1-r_0)rm = (2r_0-1)\frac{rm}{2} = (2r_0-1)\vert E\vert$.  As $r_0 > \frac{1}{2}$, $C$ has positive rate. To see that the relative distance of $C$ is positive, one uses the following result of Alon and Chung \cite[Lemma 2.3]{AC}.
  \begin{lem} In the notations of Theorem 3.2,  if $Y$ is a subset of the vertices of $X$ of size $\gamma m$, where $m=\vert X\vert$ and $0 < \gamma < 1$. Then
$$
\vert e(Y) - \frac{1}{2} r\gamma^2 m\vert ~ \leq ~\frac{r}{2} ~ \lambda ~ \gamma (1-\gamma )m
$$
where $e(Y)$ denotes the number of edges of $X$ whose both end points are in $Y$.
\end{lem}

 \begin{rem}$\frac{1}{2} r\gamma^2 m $  is what one should expect ``randomly''.
\end{rem}

Assume now that $0\neq f\in C(X,C_0)$ with edge support $D$. We want to prove that $\vert D\vert$ is large. Assume the size of $D$ is  $\frac{rm}{2} (\gamma^2 + \lambda \gamma (1-\gamma ))$ for some $0 < \gamma < 1$. Then by the Lemma, $D$ touches a set of vertices $D_0$ with at least $\gamma m$ vertices. As every edge touches two vertices, it means that on the average every vertex of $D_0$ sees at most $r(\gamma + \lambda (1-\gamma ))$ edges. So one of them sees at most this size. But it sees a vector in $C_0$ whose support is at least $\delta_0r$.  Hence $\gamma + \lambda (1-\gamma )\geq \delta_0$ which implies $\gamma\geq \frac{\delta_0-\lambda}{1-\lambda}$.  Substituting into $\vert D\vert$ implies the Theorem.

In \cite{KaW} a version of Theorem 3.1 was shown for the case when the graph $X$ is a Cayley graph of a group, in which case one can get a ``symmetric code".  This has been used in \cite{KaL} to present {\bf highly symmetric} LDPC good codes.
These codes satisfy all the ``gold standards" of coding theory: they have linear dimension (i.e., $r(c) \ge \vare > 0$), linear distance (i.e., $\delta (c) \ge \vare > 0$), they are LDPC and there exist a group $H$ acting transitively on the coordinates of $\bbf^n_2$ (i.e., acting on the Cayley graph which is an edge transitive) such that the code $C$ is invariant.  Moreover, all the constraints ($\equiv$ equations) defining $C$ are spanned by the orbit of one equation and this equation is of bounded $(\le r )$ support.

The key step for the construction of these highly symmetric codes are the edge transitive Ramanujan graphs constructed in (\cite{LSV2}) as a special case of Ramanujan complexes (\cite{LSV1}).

\section{The product replacement algorithm}

The last three decades have brought a great interest in computational
group theory. This is usually divided in two directions: one is combinatorial
group theory which deals usually with infinite groups. We will touch
this direction briefly in \S 5.2. Here we mainly deal with the
other direction: algorithms dealing with finite groups such as permutation
groups or groups of matrices over finite fields. A typical problem
in this theory is of the following type: devise an algorithm that
when given few explicit permutations in $\mathrm{Sym}\left(n\right)$
(or matrices in $\mathrm{GL}_{n}\left(\mathbb{F}_{q}\right)$) will
find various properties of the group $G$ generated by these elements,
such as:  its order, its composition factors, etc. The computational
theory of permutation groups is very developed where most problems
have deterministic algorithms. On the other hand for matrix groups
many of the practical algorithms are probabilistic.

Probabilistic algorithms very often need (pseudo) random elements from
the group $G$. Let us formulate this more formally. We need an algorithm
that when explicit elements $g_{1},\ldots,g_{k}$ (from a bigger group
like $\mathrm{Sym}\left(n\right)$ or $\mathrm{GL}_{n}\left(\mathbb{F}_{q}\right)$
) are given, it will provide us with a {}``pseudo random'' element
from $G=\left\langle g_{1},\ldots,g_{k}\right\rangle $ - the group
generated by $g_{1},\ldots,g_{k}$.

One such algorithm is to take a random word of some length $\ell$
in the generators $g_{1},\ldots,g_{k}$ and their inverses. This can
be visualized as the random walk on the Cayley graph $Cay\left(G;\left\{ g_{1}^{\pm1},\ldots,g_{k}^{\pm1}\right\} \right)$
when one stops after $\ell$ moves. This algorithm is a pretty good
one if this Cayley graph is an expander, but this is not the case
in general. The reader may think about the case $k=1$ in which case $G$
is cyclic to see how slow is the algorithm in this case.

A different approach was suggested in 1995 in \cite{CLMNO} and very quickly
became the standard way to generate random elements in finite groups
in the various packages dealing with group computations like MAGMA,
GAP etc. It is called the \emph{product replacement algorithm}.
The easiest way to describe it is also as a random walk on a graph.
This time the vertex set of the graph is $\Omega_{r}\left(G\right)=\left\{ \left(h_{1},\ldots,h_{r}\right)\in G^{r}\,\middle|\, G=\left\langle h_{1},\ldots,h_{r}\right\rangle \right\} $,
i.e., the $r$-tuples of generators of $G$.The edges correspond to
the following {}``moves'':

for $1\leq i\neq j\leq r$: \begin{eqnarray*}
L_{ij}^{\pm} & : & \left(h_{1},\ldots,h_{i},\ldots,h_{j},\ldots,h_{r}\right)\mapsto\left(h_{1},\ldots,h_{i},\ldots,h_{i}^{\pm1}h_{j},\ldots,h_{r}\right)\\
R_{ij}^{\pm} & : & \left(h_{1},\ldots,h_{i},\ldots,h_{j},\ldots,h_{r}\right)\mapsto\left(h_{1},\ldots,h_{i},\ldots,h_{j}h_{i}^{\pm1},\ldots,h_{r}\right)\end{eqnarray*}
This makes $\Omega_{r}\left(G\right)$ into a $4r\left(r-1\right)$-regular
graph. The algorithm is to take  $r > k$ and a random walk, starting at $(g_1,\ldots, g_k, e, e, \ldots, e)$,
  of say, $\ell$ steps,
then stop at a vertex of $\Omega_{r}\left(G\right)$ and pick up one
of its coordinates randomly among the $r$ possibilities. Unlike the
previous Cayley graph, this graph is highly non-symmetric, and contains
many loops and double edges. The analysis of the algorithm is very
complicated but many simulations showed outstanding performances.
For example for $G=\mathrm{Sym}\left(n\right)$, $\tau=\left(1,2\right)$
and $\sigma=\left(1,\ldots,n\right)$, the first algorithm needs (by
theoretical and experimental data) approximately $n^{2}\log n$ steps
so for $n=52$ this is over 10,000. At the same time simulations
with the product replacement algorithm for $n=52$ and $r=10$ showed
that after approximately 160 steps one gets a random-like permutation.

What is needed is a theoretical explanation for these outstanding
performances. First steps in this analysis were taken in \cite{DSC}.
A more comprehensive explanation was suggested in \cite{LP}. Here is
the crucial observation: think of $L_{ij}^{\pm}$ and $R_{ij}^{\pm}$
above as acting on the vector $\left(x_{1},\ldots,x_{i},\ldots,x_{j},\ldots,x_{r}\right)$
of $r$ free generators of the free group $F_{r}$ on $\left\{ x_{1},\ldots,x_{r}\right\} $.
Let $A^{+}=\mathrm{Aut}^{+}\left(F_{r}\right)$ be the subgroup of
$A=\mathrm{Aut}\left(F_{r}\right)$ generated by these elements. By
some well known results, going back to Nielsen, $A^{+}$ is a subgroup
of index $2$ in $A$. Now, $\Omega_{r}\left(G\right)$ can be identified
with the set $\mathrm{Epi}\left(F_{r},G\right)$ of epimorphisms from
$F_{r}$ onto $G$, where such an epimorphism $\varphi$ corresponds
to $\left(\varphi\left(x_{1}\right),\ldots,\varphi\left(x_{r}\right)\right)$.
The group $\mathrm{Aut}\left(F_{r}\right)$ acts on $\mathrm{Epi}\left(F_{r},G\right)$
by $\alpha.\varphi=\varphi\circ\alpha^{-1}$ for $\alpha\in A$. One
can easily check now that the graph structure of $\Omega_{r}\left(G\right)$
defined above is the Schreier graph of $\mathrm{Aut}\left(F_{r}\right)$
acting on the set $\Omega_{r}\left(G\right)$ w.r.t. the generators
$\big\{ L_{ij}^{\pm},R_{ij}^{\pm}\big\}$. (A Schreier graph of
a group $H$ generated by $\Sigma$ and acting on a set $X$ is the
graph with vertex set $X$ where $x\in X$ is connected to $\sigma.x$
for $\sigma\in\Sigma\cup\Sigma^{-1}$).

If the group $A=\mathrm{Aut}\left(F_{r}\right)$ has Kazhdan property
$\left(T\right)$, then an argument similar to Proposition 1.11 would give that $\Omega_{r}\left(G\right)$ are expanders.
The random walk on them converges then very fast to the uniform distribution.
This would give a conceptual explanation for the great performances
of the algorithm.

Unfortunately, it is still a (quite well known) open problem whether
$\mathrm{Aut}\left(F_{r}\right)$ has $\left(T\right)$ (it does \emph{not}
for $r=2,3$ - see \cite{GL}). Still the approach presented here was
sufficient to get some unconditional results for various classes of
finite groups:  For example for abelian groups, or more generally, nilpotent groups of bounded class.  It is shown in \cite{LP} that the subgroup of $Aut(F_r(c))$, the automorphism group of the free nilpotent group on $r$ generators  and class $c$, generated by the ``Nielsen moves" (as above) has $(T)$ if $r\ge 3$.  One can therefore deduce linear mixing time for the random walk on $\Omega_r(G)$ for $G$ nilpotent (to be compared with the subexponential bound obtained in \cite{DSC} without the use of expanders). This explains, at least for these groups, the outstanding performances of the product replacement algorithm. See \cite{LP} for details and a more general conjecture.


\chapter{Expanders in number theory}
As was mentioned (though briefly, for a more comprehensive treatment see \cite{L1}, \cite{S1}) the theory of expanders has been related to number theory in several ways. But, traditionally, the direction was from number theory to graph theory: various deep results in number theory and the theory of automorphic forms have been used to give explicit constructions of expanders and of Ramanujan graphs. We now start to see applications in the opposite direction: from expander graphs to number theory. The most notable one is the development of the affine sieve  method. This chapter will be devoted to its description and applications. For other applications see \cite{Ko1}, \cite{EHK} and \cite{EMV}.

\section{Primes on orbits}\label{sec: primes on orbits}
Many results and problems in number theory are about the existence of primes. There are infinitely many primes in $\bZ$, but Dirichlet's classical result says more:
\begin{theorem}
If $b,q\in\bZ$ with $(b,q)=1$, then there are infinitely many primes on the arithmetic progression $b+q \bZ$.
\end{theorem}
If one wants to avoid the ``local assumption" $(b,q)=1$, the result can be restated as: For every $b$ and $q\neq 0$ in $\bZ$, the arithmetic sequence $b+q \bZ$ has infinitely many numbers $x$ with $\nu(x)\leq 1+ \nu((b,q))$. Here for $x \in \bZ$ we write $\nu(x)$ for the number of prime factors of $x$.

Another well known problem about primes is:
\begin{conj}[Twin Prime Conjecture] \label{cnj:TPC}
There are infinitely many primes $p$, for which $p+2$ is also a prime.
\end{conj}

Another way to state the conjecture is: there are infinitely many $x \in \bZ$ with $\nu(x(x+2))\leq 2$.

One also expects that the Twin Prime Conjecture is true along arithmetic progressions satisfying the ``local condition" above.

A far-reaching generalization is the next conjecture of Schinzel [SS], which needs some notations: Let $\Lambda$ be an infinite subgroup of $\bZ$, i.e., $\Lambda=q \bZ$ for some $q\neq 0$, and $b \in \bZ$. Let $\cO$ be the orbit of $b$ under the action of $\Lambda$ on $\bZ$, i.e., $\cO = b+q \bZ$. Let $f(x) \in \bQ[x] $ be a polynomial which is integral on $\cO$. We say that the pair $(\cO, f)$ is \emph{primitive} if for every $2\leq k \in \bZ$ there exists $x \in \cO$ such that $(f(x),k)=1$.

\begin{conj}[Schinzel]
If $f(x) \in \bQ[x] $ is a product of $t$ irreducible factors in $\bQ[x]$, $\cO = b+q \bZ$ as above, $f$ is integral on $\cO$ and $(\cO,f)$ is primitive, then there are infinitely many $x \in \cO$ with $\nu (x) \leq t$.

\end{conj}

Taking $f(x)=x$, one gets Dirichlet Theorem and for $f(x)=x(x+2)$ the Twin Prime Conjecture in its generalized form (also along arithmetic progressions).

There are various high-dimensional conjectures generalizing Dirichlet theorem. Let us set some more notations:

Let $\Lambda$ be a non-trivial subgroup of $\bZ^n$, $b \in \bZ^n$ and $f(x_1, \dots, x_n) \in \bQ[x_1, \dots, x_n]$ which is integral on $\cO=b+\Lambda$. For $r \in \bN$, we denote by $\cO(f,r)$ the set of $x \in \cO$ for which $\nu(f(x)) \leq r$. We say that $(\cO,f)$ \emph{saturates} if for some $r<\infty$, $\cO(f,r)$ is Zariski dense in  (the Zariski closure of) $\cO$. The smallest such $r$, if it exists at all, will be denoted $r_0(\cO,f)$.

\begin{conj}[Hardy-Littlewood \cite{HL}]\label{cnj:HL}
Let $\Lambda$ be a subgroup of $\bZ^n$. Assume that for each $j$, the $j$-th coordinate function $x_j$ is nonconstant when restricted to $\Lambda$. Let $b \in \bZ^n$, $\cO = b+ \Lambda$ and $f(x)=x_1\cdot x_2\cdot \ldots \cdot x_n$, and assume $(\cO,f)$ is primitive. Then $r_0(\cO,f)=n$, i.e., the set of $x \in \cO$ whose all coordinates are simultaneously primes is Zariski dense in $b+\bC \Lambda$ and  in particular, it is infinite.
\end{conj}

       A recent breakthrough of Green, Tao and Ziegler (\cite{GTZ}, \cite{GT2})  has proved this conjecture for the case when rank$(\Lambda) \ge 2$. The most difficult case is when rank$(\Lambda) = 1$.      For example, by looking at $b=(1,3) \in \bZ^2$ and $\Lambda = \bZ(1,1)$, we see that the Twin Prime Conjecture is a special case.

Another special case is the following famous result proved not long ago by  Green and Tao \cite{GT1}:

\begin{theorem}[Arithmetic progressions of primes]\label{thm:arith prog}
For every $3\leq k \in \bN$, the set of primes contain an arithmetic progression of length $k$.
\end{theorem}

To see that Theorem \ref{thm:arith prog} is a special case of Conjecture \ref{cnj:HL}, look at $\bZ^k$ and let $\Lambda$ be the 2-dimensional subgroup $\Lambda=\bZ(1,1,1,\dots,1)+\bZ(0,1,2,3,\dots,k-1)$. Then the orbit of $(1,1,\dots,1)$ is $\Lambda$ which is the set $\{(m, m+n, m+2n, \dots, m+(k-1)n|m,n \in \bZ\}$. Conjecture \ref{cnj:HL} implies that there are infinitely many vectors of this kind whose entries are all primes.

The formulation of Hardy-Littlewood Conjecture \ref{cnj:HL}, naturally suggests to study the existence of primes vectors (i.e., vectors whose all coordinates are primes) in the orbit $\Lambda.b$ when this time $\Lambda$ is a subgroup of $GL_n(\bZ)$. Somewhat surprisingly this has not been studied till recent years. It seems that counter examples of the following kind led to think that no real theory can be developed:

\begin{exa}
Let $\Lambda$ be the cyclic subgroup of $SL_2(\bZ)$ generated by
$\left( \begin{smallmatrix}
7 & 6 \\
8 & 7\end{smallmatrix} \right)$
and $b=(1,1)^t$. The orbit $\Lambda . b$ is contained in $\{(x,y) \in \bZ^2|4x^2-3y^2=1\}$ from which one easily sees that no such $y$ is a prime, in spite of the fact that for this problem there are no ``local obstructions".
\end{exa}
Another example of a similar flavor is:

\begin{exa}[\cite{S7}, \cite{BGS2}]
Let $\Lambda= \left<\left( \begin{smallmatrix}
3 & -1 \\
1 & 0\end{smallmatrix} \right)\right> \leq SL_2(\bZ)$, $b=(2,1)^t$ and $\cO$ be the orbit $\Lambda.b$. The orbit lies on the hyperbola $\{(x,y) \in \bZ^2|x^2-3xy+y^2=1\}$ and for $n \in \bN$ we get the pairs $(f_{2n-2}, f_{2n})$, where $f_n$ is the $n$-th Fibonacci number (one can define them for $n<0$ as well). While it is conjectured that infinitely many of the $f_n$'s are primes, $f_{2n}$ are not. In fact, $f_{2n} =f_n l_n$ when $l_n$ is the $n$-th Lucas number. Moreover, it is even expected that $f_{2n}$ has an unbounded number of prime factors, when $n\rightarrow\infty$. (See \cite{BLMS}).
\end{exa}

The exciting fact, which came out only in recent years, is that these examples are the exceptional, not the typical. The Zariski closure in these cases is a torus. We will see below how conjectures and results (!) like the Hardy-Littlewood conjecture have non-abelian analogues when a torus is not involved. The key new ingredient for all this are the expanders combined with the classical combinatorial sieve of Brun. This will be our topic in the next section.

\section{Brun sieve and expanders}
For $0<x \in \bR$, denote by $\bP(x)$ the set of primes smaller or equal $x$, $\pi(x)=|\bP(x)|$ and $P(x)=\prod_{p\in \bP(x)} p$. Evaluating $\pi(x)$ is one of the most important problems in mathematics, if not the most important. Well, the prime number theorem says that $\pi(x)\sim \frac{x}{logx}$ and the Riemann hypothesis gives a sharp bound for the error term in this asymptotic result. In fact, one has an exact formula for $\pi(x)$ which was given  by Legendre at the end of the 18th century.

\begin {prop}\label{prop:8}
$$\pi(x)-\pi(\sqrt x)=-1+\sum_{S\subseteq \bP(\sqrt x )}(-1)^{|S|}\Big \lfloor \frac{x}{\prod_{p\in S}p}\Big\rfloor $$
\end{prop}
By nowadays's  standards the proof is a simple application of the inclusion-exclusion formula: the primes between $\sqrt{x}$ and  $x$ are those integers $n\neq 1$ which are not divisible by any prime less than $\sqrt{x}$. So we count them by taking $x$, subtracting those divisible by one prime less than $\sqrt{x}$, adding the number of those divisible by two primes, etc. Basically, we are applying the classical Eratosthenes' sieve method.

While Proposition \ref{prop:8} gives an exact formula, it is not very useful. The error term, for example, between $\frac{x}{\prod_{p\in S}p}$ and its integral part is ``small"-bounded by $1$. But there are so many summands (approximately $4^{\sqrt{x}/log x}$) which makes the formula quite useless.

Various ``sieve methods" have been developed for problems like that--the reader is referred to \cite{FI} and \cite{IK}, for example. Let us say a few words about the combinatorial sieve developed by Brun. His main motivation was to handle the Twin Prime Conjecture \ref{cnj:TPC}. In a different language it says that if $f(x)=x(x+2)$ then for infinitely many $n$'s in $\bN$, $f(n)$ is a product of only two primes, i.e. $\nu(f(n))\leq 2$.
Let $f(x)$ be any integral polynomial $f(x) \in \bZ[x]$ e.g., $f(x)=x(x+2)$. Let $x$ be a large real number, and $z<x$. Denote:
$$S(f,z):=\mathop{\sum_{n\leq x}}\limits_{(f(n),P(z))=1}1$$
so $S(f,z)$ counts those $n$ less than $x$, such that $f(n)$ is not divisible by any prime less than $z$. Of course, we want $z$ to be as large as possible.

Recall the Mobius function
\begin{equation*}
\mu(n) = \left\{
\begin{array}{rl}
1 & \text{if } n=1\\
(-1)^r & \text{if } n=p_1\cdot\ldots\cdot p_r \text{  distinct primes }\\
0 & \text{otherwise}
\end{array} \right.
\end{equation*}
The following is well-known and easy to prove:
\begin{equation*}
\sum_{d|n}{\mu(d)} = \left\{
\begin{array}{rl}
1 &  n=1\\
0 & n>1
\end{array} \right.
\end{equation*}
One can therefore now write:
$$S(f,z):=\mathop{\sum_{n\leq x}}\limits_{(f(n),P(z))=1}1=\sum_{n\leq x}\sum_{d|(f(n),P(z))}\mu(d)=$$ $$=\sum_{d|P(z)}\mu(d)\left(\mathop{\sum_{n\leq x}}\limits_{f(n)\equiv 0(d)}1\right)$$
We now denote:
$$\beta(d)=|\{m \text{ mod } d | f(m)\equiv0(d)\}|$$
i.e., the  number of solutions of $f$ mod $d$. Running over all $n$'s up to $x$, then  one runs approximately $\frac{x}{d}$ times on all the residues $mod$ $d$, approximately $\frac{x}{d} \beta(d)$  of them will give zeros for $f$ mod $d$. So, continuing the evaluation of $S(f,z)$ we have:
$$S(f,z)=\sum_{d|P(z)}\mu(d)(\frac{\beta(d)}{d}x+r(d))$$
where $r(d)$ is an error term. Note that $\frac{\beta(d)}{d}$ is a multiplicative function of $d$.

Brun developed a method to analyze such a sum with  particular interest in the case $f(x)=x(x+2)$. He deduced that for $\delta$ small enough if $z=x^\delta$ then $S(f,z)\geq c \frac{x}{log(x)^2}$ which means that there are infinitely many $n$'s with no prime divisor for $f(n)$ less than $n^\delta$. For such $n$'s, $\nu(f(n))\leq \frac{\deg f}{\delta}$. So, while he felt short from proving the Twin Prime Conjecture he was able to show that there are infinitely many $n$'s with $\nu(f(n))\leq 18$. His method has been refined; the current record is due to Chen \cite{Ch} who replaced 18 by 3. Namely, there are infinitely many pairs $(n,n+2)$ such that one of them is a prime and the other is a product of at most two primes.

These ``combinatorial sieve methods" have also been applied to the higher dimensional cases described in \S \ref{sec: primes on orbits}. For examples, one gets a partial result toward Hardy-Littlewood Conjecture: in the notation of Conjecture \ref{cnj:HL}, on can prove that there exists $r\in \bN$ such that $r_0(\cO,f)\leq r$. In particular, the orbit $b+\Lambda$ contains infinitely many vectors all whose entries are product of a bounded number of primes. Moreover, it has even been proved that this $r$ depends only on $n$ and not on $b$ or $\Lambda$ (assuming, of course, no local obstructions, i.e. $(\cO,f)$ is primitive).

All this is a quite deep (and quite technical) theory. The relevant for our story came from the insight of Sarnak who noticed that the machinery of the Brun sieve can be carried out also for a general non-commutative subgroups $\Lambda \leq GL_n(\bZ)$ acting of $\bZ^n$, \emph{provided} $\Lambda$ has property ($\tau$) w.r.t congruence subgroups. The orbit in this case of $b\in \bZ^n$ is $\Lambda .b=\{\gamma .b|\gamma\in \Gamma\}$, and one can start the same kind of computation we illustrated above for the Twin Prime problem. This time, instead of summing up over all $n\leq x$, one sums over the ball of radius at most $k$, with respect to a fixed set of generators $\Sigma$ of $\Lambda$. The crucial point is that these balls $B(k)=\{\gamma\in \Lambda| \text{length}_\Sigma(\gamma)\leq k\}$ when act on $b\in \bZ^n$ and reduced mod $d$, i.e., the set $B(k).b(\text{mod } d)$, distribute almost uniformly over the vectors $\Gamma. b (\text{mod } d)$, as a  subset of $(\bZ/d\bZ)^n$. This is exactly what the expander property gives us (compare with Proposition 1.6).

At first sight this connection with expanders may look counter intuitive: we want to extend sieve methods from abelian cases, such as Hardy-Littlewood Conjecture, to a non-abelian setting. The abelian case \emph{never} gives rise to expanders (see \cite{LW}) - why should this be the needed property in the non-abelian case? The point is that in the abelian setting the number of integer points in arithmetic progressions which are contained in a large interval can be estimated quite accurately in the obvious way. But in the non-abelian setting, it is not clear what is the distribution of the points in a ball when taken mod $d$. Note also that such groups have usually exponential growth and so when we move from ball $B(k)$ to $B(k+1)$ the boundary is as large as the original ball. The expanding property enables one to overcome this difficulty. In fact, one does not need that $\Lambda < GL_n(\bZ)$ has ($\tau$) with respect to all congruence subgroups, it suffices to know it with regard to congruence subgroups mod $d$ when $d$ is square-free.

All this machinery was put to work in the paper of Bourgain-Gamburd-Sarnak \cite{BGS2}. At the time when the paper was written property $(\tau)$ was known for such $\Lambda$'s only when the Zariski closure of $\Gamma$ was $SL_2$ (due to Helfgott \cite{H}, Bourgain-Gamburd \cite{BG1} and the extension to all square-free in \cite{BGS2}). But they also proved some conditional results, assuming an affirmative answer to some generalized form of the 1-2-3 problem (See \S 2.5). That work gave a push to efforts in this direction by a good number of authors (\cite{BG3},\cite{BG4}, \cite{BGS3}, \cite{BV}, \cite{V}, \cite{BGT2}, \cite{PS2},  \cite{S7}, \cite{SGV}). The most general result for the ``affine-sieve method" as it is called now, is given in a forthcoming paper of Salehi-Golsefidy and Sarnak \cite{SGS}:

\begin{theorem}\label{thm:general}
Let $\Lambda \leq GL_n(\bZ)$ be a finitely generated subgroup with Zariski closure $G$. Assume the reductive part of $G^0$ -- the connected component of $G$ -- is semisimple. Let $b\in \bZ^n$, $\cO=\Gamma .b$ and $f\in \bQ[x_1,\dots,x_n]$ which is integral and not constant on $\cO$. Then $(\cO,f)$-saturates, namely there exist $r\in \bN$ such that the set of vectors in $\cO$ whose components are product of at most $r$ primes is Zariski dense in $\cO$.
\end{theorem}

This theorem covers also cases when $G$ is unipotent (and so various classical results) as well as completely new cases when $G$ is semisimple. The method is called ``affine sieve" as it also covers ``affine transformations" of $\bZ^n$ and not only linear. The affine case can be easily reduced to a linear case of higher dimensions. Some of the classical problems (e.g. the Hardy-Littlewood Conjecture) are naturally expressed as affine problems rather than linear.

The case which is not covered by the last general theorem is of a torus (or when one has central torus in $G$). Some of the difficult problems in number theory can be presented in this language: e.g. the Mersenne Conjecture: there are infinitely many primes $p$ with $2^p -1$ also a prime, is such a problem (See \cite[\S 2.1]{BGS2}). But the set of primes is ``too thin" to sieve over it. So the new method shed no new light on this conjecture. It is not even known if there are infinitely many almost primes of the form $2^n-1$.

Still there are few concrete problems where the new method gives some fascinating results. Some of them will be described in the next section.

\section{Some applications to classical problems}\label{sec:applications}
Theorem \ref{thm:general} above gives some results which are completely out of reach by other methods. E.g., if $\Lambda\leq GL_n(\bZ)$ is a group as in the theorem, then the group itself contained infinitely many matrices which are almost primes, i.e., all their entries are products of a bounded number of primes. An example satisfying this is \emph{any} non-virtually cyclic subgroup of $SL_2(\bZ)$. But, here Theorem \ref{thm:general} answers questions which have not asked before.

Let us now present (following \cite{BGS2}, \cite{S5}, \cite{S8}) two applications to classical number theoretic problems:

\subsection*{Pythagorian triangles} Look at right angle triangles with integral edges $x_1,x_2$ and $x_3$ so $ x_3^2=x_1^2+x_2^2$ and assume $g.c.d(x_1,x_2,x_3)=1$. It is well known that in this case there exist $m,n\in \bZ$, one odd, one even and $(m,n)=1$ s.t. $x_1=m^2-n^2$,$x_2=2mn$ and $x_3=m^2+n^2$. It follows that $x_1$ is divisible by $3$ and $x_2$ by $4$. So the area of the triangle $A=\frac{x_1x_2}{2}$ is divisible by $6$. All the Pythagorian triples $(x_1,x_2,x_3)$ as above are obtained as the orbit $\cO$ of the triples $(3,4,5)$ acted upon by $O_F(\bZ)$ when $F$ is the quadratic form $x_1^2+x_2^2-x_3^2$ and $\Lambda=O_F(\bZ)$ is the group of $3\times 3$ integral matrices preserving this form. The group $\Lambda$ satisfies the conditions of Theorem \ref{thm:general} and $f=\frac{x_1x_2}{2}$ is integral on $\cO$ (and even divisible by $6$). We deduce that there are infinitely many triples whose areas are almost primes.

Now what is $r_0(\cO,f)$ - i.e. what is the minimal $r$ for which the set of triples with $\nu(\text{area})\leq r$ is Zariski dense? This is a more delicate question. Some elementary arguments show that it is at least $6$ and some recent work of Green and Tao \cite{GT2} implies that it is indeed $6$. (See \cite{BGS2} and the references therein for more information).
\subsection*{Integral Apollonian Packing}
\begin{figure}[htp]
\centering
\includegraphics[scale=0.1]{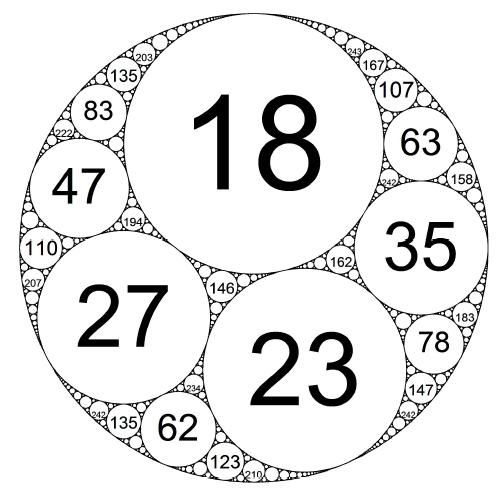}
\caption{Apollonian packing}\label{fig:ap1}
\end{figure}

A classical theorem of Apollonius asserts that given three mutually tangents circles $C_1,C_2$ and $C_3$, there are exactly two circles $C_4$ and $C_4'$ tangents to all three. Descartes' Theorem says that the curvatures of these circles (i.e. the reciprocals of the radii) $a_1,a_2,a_3,a_4$ satisfy $F(a_1,a_2,a_3,a_4)=0$ where
\begin{equation}
F(a_1,a_2,a_3,a_4)=2(a_1^2+a_2^2+a_3^2+a_4^2)-(a_1+a_2+a_3+a_4)^2
\label{eq:packing}
\end{equation}
(a negative solution correspond to a situation when one circle touches the others from the outside). An easy calculation using \eqref{eq:packing} shows that given $C_1,C_2,C_3$ with curvatures $a_1,a_2,a_3$ respectively, there are two solutions $C_4$ and $C_4'$ with curvatures $a_4$ and $a_4'$ satisfying
\begin{equation}
a_4'=2a_1+2a_2+2a_3-a_4
\label{eq:packing2}
\end{equation}
It also shows that starting with an integral vector $(a_1,a_2,a_3,a_4)$, the other quadruple $(a_1,a_2,a_3,a_4')$ is also integral. This can be carried out with any sub-triple of $C_1,C_2,C_3,C_4$. We deduce: let

\begin{align*}
S_1=&\Biggl(\begin{smallmatrix}
\!\!-1 & 2 &2 &2 \\
0 & 1 &0 &0 \\
0 & 0 &1 &0 \\
0 & 0 &0 &1 \end{smallmatrix}\Biggr),\qquad  \; S_2= \;  \Biggl(\begin{smallmatrix}
1 & 0 &0 &0 \\
2 & \!\!-1 &2 &2 \\
0 & 0 &1 &0 \\
0 & 0 &0 &1 \end{smallmatrix}\Biggr), \\
S_3= &\Biggl(\begin{smallmatrix}
 1 & 0 &0 &0\\
0 & 1 &0 &0 \\
2 & 2 \!\!&-1 &2 \\
0 & 0 &0 &1 \end{smallmatrix}\Biggr) \,  \text{ and } \, S_4= \Biggl(\begin{smallmatrix}
1 & 0 &0 &0 \\
0 & 1 &0 &0 \\
0 & 0 &1 &0 \\
2 & 2 &2 \!\!&-1 \end{smallmatrix}\Biggr),
\end{align*}
and let $\Lambda$ be the subgroup generated by these four reflections. Then starting with any integral quadruple $b=(a_1,a_2,a_3,a_4)$ of integral curvatures of mutually touching circles, all elements in the orbit $\Lambda . b$ represent such quadruples. Moreover if $\gamma'=S_i\gamma$ in $\Lambda$, then the corresponding quadruple share a common triple. See Figure \ref{fig:ap1} for the starting stages of the orbit $(18,23,27,146)$.

The subgroup $\Lambda$ of $GL_4(\bZ)$ preserves the quadratic form $F$ of equation \eqref{eq:packing}. It therefore lies within a conjugate of $SO(3,1)$ and in fact, it is Zariski dense there. One can therefore deduce from Theorem \ref{thm:general} various number theoretic results on the orbit $\Lambda . b$.

But many more questions come up naturally. Are there infinitely many primes in the set of curvatures of the circles in the orbit? How many? One wishes to have a ``prime number theorem" estimating the density of prime curvatures within the orbit of the ball of radius $N$ in $\Lambda$ (w.r.t. $\{S_i\}_{i=1}^4$). Are there infinitely many twin prime? i.e., are there infinitely many kissing pairs of circles with prime curvatures.

A rich theory has started to emerge in recent years (cf. \cite{GLMWY}, \cite{S5}, \cite{S8}, \cite{KO1}, \cite{Fu}, \cite{BF}, \cite{FS} and the reference therein). This is a fascinating crossroad of number theory, geometry, group theory, dynamics, ergodic theory and expanders!

\chapter{Applications to group theory}

In the previous chapters the connections between expander graphs and
group theory has been illustrated over and over again. Many (in some
sense {}``most'') of the examples we gave for expander graphs were
Cayley graphs of groups and being expanders says something quite deep
on their group structure and/or their representation theory.

In this chapter we will describe several results about groups whose formulations does not mention expanders but expanders come out substantially
in the proofs, sometimes in a somewhat surprising way. We will start
with describing a new sieve method for finitely generated groups and
indicate several applications to linear groups and to the mapping
class groups. This is a new direction and it can be expected that
this sieve method will have further use in group theory.

We will also bring some application to combinatorial group theory
- a direction whose significance will be fully appreciated in the
next chapter when we will discuss geometric applications.

\section{Measuring subsets of finitely generated groups}

Let $G=\mathrm{GL}_{n}\left(\mathbb{C}\right)$ the group of $n\times n$
invertible complex matrices. A standard claim on $G$ is: For generic
$g\in G$ the centralizer $C_{G}\left(g\right)$ of $g$ in $G$ is
abelian. What do we mean by this statement? What is meant by {}``generic''?
Well, one can work in the Baire category setting, in the measure theoretic
language or in the Zariski topology. Whatever setting we choose the
statement says that outside of a {}``meager'' subset of $G$, the
property of abelian centralizer is satisfied. The proof is easy: for
almost every $g\in G$ the eigenvalues of $g$ are all distinct (since
the set of zeros of the discriminant is {}``meager''). Thus, with
a suitable basis of $\mathbb{C}^{n}$, $g$ is diagonal with $n$
different eigenvalues and the centralizer of such $g$ is just the
diagonal matrices and hence abelian.

Let us now look at the finitely generated group $\Gamma=\mathrm{SL}_{n}\left(\mathbb{Z}\right)$.
If we want to claim a similar statement about $\Gamma$: {}``For
generic $g\in\Gamma$ the centralizer $C_{\Gamma}\left(g\right)$
is abelian''. What does this mean? Is there a natural way {}``to
measure'' subsets of $\Gamma$ to be {}``small'' or {}``large''?
On a countable set like $\Gamma$ one cannot take {}``a uniform distribution''.
This is a problem not only in group theory. The reader is referred
to an interesting lecture by Barry Mazur \cite{Maz} which illustrates
via various questions in number theory, that one may have quite different
answers to similar problems depending on the probabilistic model.

These types of questions have recently received attention also from another point of view.  Complexity theory, i.e., the theory of algorithms, usually deals with ``worse case scenario", i.e. a problem is considered ``difficult" if it is difficult for \emph{some} inputs.  But, in real life quite often we really care whether it is easy or difficult for ``most" cases, for the ``generic" inputs.  These two different approaches to complexity theory can be very different. There are even ``undecidable" problems which can be solved in polynomial time for ``most" inputs.  In recent years there have been a number of papers studying this direction in combinatorial group theory.  Some of it was motivated by the proposal of various crypto systems based on the Braid groups.

This led to various approaches to the notion of ``generic" elements in a finitely generated group (cf. \cite{KMSS1}, \cite{KMSS2}, \cite{KS1} and \cite{MR}). Let us call the reader's attention  to \cite{BHKLS}, where it is shown that the answer to a problem can be very different in two different models of randomness, even if both are ``natural".

Anyway, here is the model we will work with: Let $\Gamma$ be a group
generated by a finite symmetric set $\Sigma$. We will assume that
$\Sigma$ satisfy some relation of odd length (This is a non-essential
condition that simplify the notations, avoiding bipartite graphs in
what follows. It happens automatically if $e$ - the identity element
of $\Gamma$ - is in $\Sigma$). A walk $w$ on $\Gamma$ w.r.t. $\Sigma$
is a function $w:\mathbb{N}^{+}\rightarrow\Sigma$. The $k$-th step
of $w$ is $w_{k}:=w\left(1\right)\cdot\ldots\cdot w\left(k\right)$
(where $w_{0}=e$). The uniform measure $\mu$ on $\Sigma$ induces
a product measure $\overline{\mu}$ on the set of $\Sigma$-walks
$W_{\Sigma}:=\Sigma^{\mathbb{N}^{+}}$. For a subset $Z$ of $\Gamma$
we denote the probability that the $k$-th step of a walk belongs
to $Z$ by $prob\left(w_{k}\in Z\right)$. We say that $Z$ is \emph{exponentially
small w.r.t. $\Sigma$} if there are constants $c,\alpha>0$ s.t.
$prob\left(w_{k}\in Z\right)\leq ce^{-\alpha k}$ for every $k\in\mathbb{N}$.
If this happens w.r.t. every such $\Sigma$  (where $c$
and $\alpha$ may depend on $\Sigma$), $Z$ is \emph{exponentially
small}. We will say that $Z$ is \emph{exponentially generic}
if its complement in $\Gamma$ is exponentially small.

Now, once there
is a meaning to be small and large we can see that the set of $g$'s
in $\Gamma=\mathrm{SL}_{n}\left(\mathbb{Z}\right)$ for which $C_{\Gamma}\left(g\right)$
is not abelian is exponentially small. Indeed, fix a set $\Sigma$
of generators for $\Gamma$ and let $k\in\mathbb{N}$ and $p$ be
a prime of size exponential in $k$. The set $Z_{p}$ of matrices
in $\mathrm{SL}_{n}\left(\mathbb{F}_{p}\right)$ with multiple eigenvalues
(i.e. discriminant equal zero) satisfies $\frac{\left|Z_{p}\right|}{\left|\mathrm{SL}_{n}\left(\mathbb{F}_{p}\right)\right|}\sim\frac{1}{p}$.
By Proposition 1.12 and Theorem 2.3, $Y=Cay\left(\mathrm{SL}_{n}\left(\mathbb{F}_{p}\right);\Sigma\right)$
are $\varepsilon$-expanders for some $\varepsilon$ depending on
$\Sigma$ but not on $p$. Thus by Proposition 1.6, the random walk
on $Y$ falls into $Z_{p}$ at time $k$ with probability approximately $\frac 1p$  which is
exponentially small in $k$.

The above argument illustrates how expanders play an important role
in measuring subsets of $\Gamma$.

Similar ideas can lead to results which on first sight look very different:

\medskip{}

\begin{theorem}[\cite{BCLM}] Let $\Gamma$ be a linear group
generated by a finite set $\Sigma$. If $\Gamma$ is not virtually
nilpotent then $\Gamma$ has exponential conjugacy growth, i.e., there
exists a constant $C>1$ such that for every $k\gg0$, the ball of
radius $k$ around the identity in $Cay\left(\Gamma;\Sigma\right)$
intersects non-trivially at least $C^{k}$ different conjugacy classes.
\end{theorem}

Here the point of the proof is that in congruence quotients, each
conjugacy class is {}``small''. In the last result, and in the stronger
forms of it in \cite{BCLM}, one does not really need the full power
of expanders and results like Theorem 2.22 above suffice. That is why
it holds also in positive characteristic. This is not the case for
the more powerful method we will apply in the next section which need
also Theorem 2.24. At this point, this last result is known only in
zero characteristic. It will be very useful to prove an analogous
result for positive characteristic.

\section{Sieve method in group theory}

The results mentioned in the previous section {}``measure'' a subset
$Z$ of $\Gamma=\mathrm{SL}_{n}\left(\mathbb{Z}\right)$ (or more
general linear groups) by projecting it to $\mathrm{SL}_{n}\left(\mathbb{F}_{p}\right)$
for one prime $p$, at a time, showing that the projection $Z_{p}$
is {}``small'' and since $\mathrm{SL}_{n}\left(\mathbb{F}_{p}\right)$
is an expander the random walk meets $Z_{p}$ with exponentially small
probability. This method works for sets $Z$ for which the projections
$Z_{p}$ are small, e.g., when $Z$ is an intersection of an algebraic
variety $V$ with $\Gamma$. In this case, indeed, the projection
$Z_{p}$ is {}``small'' by the famous Lang-Weil Theorem. But for
various natural problems the projection $Z_{p}$, of $Z$ mod $p$
is large (say, proportional to the size of $\mathrm{SL}_{n}\left(\mathbb{F}_{p}\right)$).
For dealing with such problems one needs the group sieve method, which
is a group theoretical analogue of the {}``large sieve'' in analytic
number theory.

We can now formulate the general {}``group sieve method'':

\begin{theorem}\label{sieve theorem} Let $\Gamma$ be a finitely generated group. Let $(N_i)_{i \in I}$ be a series of
finite index normal subgroups of $\Gamma$ where $I \subseteq \bbn$.
Assume that there are constants $\gamma>0$ and $d \in \bbn^+$ such
that:
\begin{itemize}
\item[1.] $|\{i \in I \mid i \le e^{k}\}|\ge e^{\gamma k}$ for every large enough $k \in
\bbn$.
\item[2.] $\Gamma$ has property $(\tau)$ w.r.t. the family of
normal subgroups $(N_i \cap N_j)_{i,j\in I}$.
\item[3.] $|\Gamma_i| \le i^d$ for every $i \in I$ where $\Gamma_i:=\Gamma/N_i$.
\item[4.] The natural map $\Gamma_{i,j}\rightarrow \Gamma_i \times
\Gamma_j$ is an isomorphism   for every distinct $i,j\in I$ where
$\Gamma_{i,j}:=\Gamma /N_i \cap N_j$.
\end{itemize} Then a subset $Z\subseteq \Gamma$ is exponentially small if there is $c>0$ such that:
\begin{itemize}
\item[5.] $\frac{|Z_i|}{|\Gamma_i|}\le 1-c$ for every $i \in I$ where $Z_i:=ZN_i/N_i$.
\end{itemize}
\end{theorem}

The above formulation is taken from \cite{LM1}. This is a generalization
(and simplification) of a method used by Rivin \cite{Ri},
and by Kowalski \cite{Ko1}. It has been greatly influenced by the
{}``affine sieve'' of Chapter 4.

The last theorem gives a very general result but applying it for particular
cases still requires a substantial amount of work. The more difficult
part is establishing properties 2 and 5 in the theorem. Property
2 is true for a large class of groups by Theorem 2.24. (This theorem
is formulated for subgroups of $\mathrm{GL}_{n}\left(\mathbb{Q}\right)$
but usually one can reduce questions about general finitely generated
linear groups, over fields of characteristic $0$, to this case by the
method of specialization cf. \cite{LM1}). As mentioned in \S 5.1, it will
be useful to have an analogue of Theorem 2.24 for fields of positive
characteristic. Once this will be done the group sieve method should
give various applications also for these cases.

Property 5 of Theorem 5.2 depends very much on each specific problem.
Let us describe here the case where $Z\subseteq\Gamma$ is the subset
of all proper powers in $\Gamma$, i.e. $Z=\bigcup_{m\geq2}Z\left(m\right)$
where for $2\leq m\in\mathbb{N}$, we denote $Z\left(m\right)=\left\{ \gamma^{m}\,\middle|\,\gamma\in\Gamma\right\} $
the \emph{set} of $m$-powers. The main result of \cite{LM1} is
the following:

\begin{theorem} Let $\Gamma$ be a finitely generated linear
group over a characteristic 0 field. Assume $\Gamma$ is not virtually
solvable. Then \be*Z=\bigcup\limits _{2\leq m\in\mathbb{N}}Z\left(m\right)=\bigcup\limits _{2\leq m\in\mathbb{N}}\left\{ \gamma^{m}\,\middle|\,\gamma\in\Gamma\right\} \ee*
is an exponentially small subset of $\Gamma$.
\end{theorem}

This theorem is a far reaching straightening of the main result of
\cite{HKLS}. Not only it gives a quantitative result on $Z$, but it
also deals with the union of all the $Z\left(m\right)$'s together. In
\cite{HKLS} only finitely many $m$'s could be considered at a time.
This is the power of the sieve which enables such a stronger result.

The proof of property 5 of Theorem 5.2 for this case also needs some
careful treatment: the projection of $Z$ is \emph{onto} for \emph{every}
finite quotient. Thus one has to treat each $Z\left(m\right)$ separately,
getting quantitative results and then summing them together (see \cite{LM1}
for details).

\section{The mapping class group}

In this section we will apply the group sieve method to $A=\mathrm{Aut}\left(F_{n}\right)$
the automorphism group of the free group on $n$ generators and to
$M=\mathrm{MCG}\left(g\right)$ the mapping class group of a closed surface $S_{g}$ of genus $g$. The group $M$ is isomorphic
to $\mathrm{Out}\left(\Pi_{g}\right)={\mathrm{Aut}\left(\Pi_{g}\right)}/{\mathrm{Inn}\left(\Pi_{g}\right)}$,
the group of outer automorphisms of $\Pi_{g}=\pi_{1}\left(S_{g}\right)$
the fundamental group of $S_{g}$. The group $\Pi_{g}$ has a presentation
with $2g$ generators $a_{1},\ldots,a_{g},b_{1},\ldots,b_{g}$ subject
to one relation $\prod_{i=1}^{g}\left[a_{i},b_{i}\right]$ where $\left[a,b\right]=a^{-1}b^{-1}ab$.
The mapping class group is of great importance in topology and geometry
and we will come back to these aspects in Chapter 6 where we will
treat geometric applications of expanders. Here we mainly treat it
from its algebraic description, though, the major question comes from
topology.

Thurston classified the elements of $M$ into three kinds $\left(i\right)$
pseudo-Anosov $\left(ii\right)$ reducible and $\left(iii\right)$
elliptic. This is somewhat similar in spirit to the classification
of elements of $\mathrm{MCG}\left(1\right)=\mathrm{SL}_{2}\left(\mathbb{Z}\right)$
into hyperbolic, parabolic and elliptic. We will not give the exact
definitions sending the reader to \cite{Ri}, \cite{Ko1} and the references therein for details. He conjectured
that {}``generic'' elements of $M$ are pseudo-Anosov. In one
form (which is weaker and stronger than the following theorem)  this was proved by Maher (\cite{Ma1}, \cite{Ma2}). Rivin \cite{Ri} (see also Kowalski \cite{Ko1}) proved, by using the sieve
method:
\begin{theorem} The set of pseudo-Anosov elements of $M=\mathrm{MCG}\left(g\right)$
is exponentially generic.
\end{theorem}
The proof uses the fact that $M=\mathrm{MCG}\left(g\right)$ is mapped
onto the arithmetic group $\Gamma=\mathrm{Sp}\left(2g,\mathbb{Z}\right)$.
Now, a criterion due to Casson and Bleiler gives a sufficient condition
for an element $\gamma$ of $M$ to be pseudo-Anosov in terms of some
conditions on its image $\gamma'\in\Gamma$. Rivin showed that the
set of those $\gamma'\in\Gamma$ which do not satisfy this condition
is exponentially small, and deduced that the non pseudo-Anosov elements
of $M$ form an exponentially small set.

The proof sketched above gives no information of the subgroup $T=\break\ker\left(\mathrm{MCG}\left(g\right)\rightarrow\mathrm{Sp}\left(2g,\mathbb{Z}\right)\right)$
- the Torelli subgroup of the mapping class group. It was asked by
Kowalski \cite{Ko1} whether a similar result to Theorem 5.4 holds also
for $T$. In \cite{LM2} and in \cite{MS}, independently, it was shown to be the case by using various
representations of $T$ onto $\mathrm{Sp}\left(2\left(g-1\right),\mathbb{Z}\right)$
obtained by considering the action of $T$ on the homology of the
2-sheeted covers of the surface $S_{g}$ (or equivalently on the commutator
quotients of the index 2 subgroups of $\Pi_{g}$) in the spirit of \cite{Lo} and \cite{GL}.

The above mentioned results have analogous results, proved also in \cite{Ri}, \cite{Ko1} and \cite{LM2}, for $\mathrm{Aut}\left(F_{n}\right)$
replacing $\mathrm{MCG}\left(g\right)$.The role of pseudo-Anosov
is played by the {}``fully irreducible'' automorphisms (called also:
\emph{irreducible with irreducible powers} - iwip, for short).
These are the automorphisms $\alpha\in\mathrm{Aut}\left(F_{n}\right)$
such that no positive power of $\alpha$ sends a free factor $H$
of $F_{n}$ to a conjugate. The conclusion is that these automorphisms
are exponentially generic in $\mathrm{Aut}\left(F_{n}\right)$ as
well as in $IA\left(n\right)=\ker\left(\mathrm{Aut}\left(F_{n}\right)\rightarrow\mathrm{GL}_{n}\left(\mathbb{Z}\right)\right)$,
when $n\ge 3$.

\section{The generic Galois group of linear groups}

The method of proof of the results in the previous section showed
(when establishing the criterion of Casson and Bleiler mentioned there)
something stronger which is of independent interest: for an exponentially
generic matrix $A\in\mathrm{SL}_{n}\left(\mathbb{Z}\right)$, the
Galois group over $\mathbb{Q}$ of the splitting field of the characteristic
polynomial of $A$ is isomorphic to $\mathrm{Sym}\left(n\right)$,
the full symmetric group on $n$ letters. Similarly, for generic elements
in $\mathrm{Sp}\left(2g,\mathbb{Z}\right)$ the Galois group of the
splitting polynomial is isomorphic to the Weyl group of the algebraic
group $\mathrm{Sp}\left(2g\right)$. A common generalization was proved
by Jouve, Kowalski and Zywina \cite{JKZ}.

\begin{theorem} Let $k$ be a number field and $\mathbf{G}$
a connected semisimple group defined and split over $k$ with a faithful
representation $\rho:\mathbf{G}\rightarrow\mathrm{GL}\left(m\right)$
defined over $k$. Let $\Gamma\subseteq\mathbf{G}\left(k\right)$
be an arithmetic subgroup. Then for exponentially generic elements
$A$ in $\Gamma$, the Galois group of the splitting field over $k$
of the characteristic polynomial of $A$ is isomorphic to the Weyl
group $W\left(\mathbf{G}\right)$ of the algebraic group $\mathbf{G}$.
\end{theorem}

Although this statement seems to be asymptotic, the method is effective and enables one, for example, to find matrices whose characteristic polynomials have $W(E_8)$ as their Galois groups over $k$.

The reader is referred to \cite{JKZ} for a more general result when
$\mathbf{G}$ does not split and to \cite{LR} for more general linear
groups. (The results are somewhat different!)

Those results use heavily the fact that congruence quotients are expanders.
But they need also an interesting use of Chebotarev theorem which
{}``provides'' elements in conjugacy classes of the {}``target''
Galois group. These elements are defined only up to conjugacy. This
leads to the following interesting notion:

\begin{definition} A subset $S$ of a finite group $G$ is
said to generate $G$ invariably if $G=\left\langle s^{g\left(s\right)}\,\middle|\, s\in S\right\rangle $
for \emph{any} choice of $g\left(s\right)\in G$. (i.e. if every
element of $S$ is replaced by some conjugate of it, we still get
a set of generators).
\end{definition}

 This is an interesting group theoretic invariant of importance for computational group theory. For some basic properties of it see \cite{KLS} and the references therein. It illustrates once again how results and methods from pure mathematics and computer science enrich each other back and forth.

\section{Property $\left(\tau\right)$ in combinatorial group theory}

Let $\Gamma$ be a discrete group. It is called \emph{residually finite}
if the intersection of the finite index subgroups is trivial. We say
that $\Gamma$ \emph{splits} if $\Gamma$ can be written as free
product with amalgamation $A\underset{C}{*}B$ or as an HNN-construction
$A_{*_{C_{1}=C_{2}}}$ in a non-trivial way, i.e., $C\lneqq A,B$. It is well known that $\Gamma$
splits if and only if it acts on a simplicial tree without a (common)
fixed point. Note that if $\Gamma$ is finitely generated then it is an HNN-construction    if and only if
$\Gamma$ is mapped surjectively onto the infinite cyclic group $\mathbb{Z}$ and this happens iff the commutator subgroup $\left[\Gamma,\Gamma\right]$
of $\Gamma$ is of infinite index.

For a finitely generated group $\Lambda$ we denote by $d\left(\Lambda\right)$
- the minimal number of generators of $\Lambda$. The \emph{rank gradient}
of $\Gamma$, $RG\left(\Gamma\right)$ is defined as:\[
RG\left(\Gamma\right)=\inf\left\{ \frac{d\left(\Lambda\right)-1}{\left[\Gamma:\Lambda\right]}\,\bigg|\,\mathrm{{\Lambda\, finite\, index\atop subgroup\, of\,\Gamma}}\right\} \]
The following result of Lackenby \cite{La1} gives a surprising connection
between $\left(\tau\right)$, splitting and $RG\left(\Gamma\right)$.

\begin{theorem} Let $\Gamma$ be a finitely presented residually
finite group. Then $\Gamma$ satisfies (at least) one of the following
three properties:
\begin{description}
\item [{$\left(a\right)$}] $\Gamma$ virtually splits (i.e. has a finite
index subgroup $\Lambda$ which splits).
\item [{$\left(b\right)$}] $\Gamma$ has property $\left(\tau\right)$.
\item [{$\left(c\right)$}] $RG\left(\Gamma\right)=0$.
\end{description}
\end{theorem}

The method of proof goes like that: one assumes that $\Gamma$ does
not have $\left(b\right)$ and $\left(c\right)$, i.e. the quotient
graphs of $\Gamma$ are \emph{not} expanders and $RG\left(\Gamma\right)>0$
which means that the number of generators of finite index subgroups
grows linearly with the index. These two pieces of information are
used to deduce that a suitable finite cover $Y$ of a 2-dimensional
complex $X$ with $\pi_{1}\left(X\right)=\Gamma$ can be decomposed
 as $Y=Y_{1}\cup Y_{2}$ in a non-trivial way that will enable to apply van Kampen
Theorem to deduce that $\pi_{1}\left(Y\right)$ splits - see \cite{La1}
for details.

As we will see in the next chapter, it is of great importance in the
theory of 3-manifolds to be able to show that $\pi_{1}\left(M\right)$
of such a manifold $M$, virtually splits. So a result like Theorem
5.7 and various variants of it, are useful there as a tool to get
the desired conclusion.

Another application is that for every finitely presented amenable
group $\Gamma$, $RG\left(\Gamma\right)=0$ since such $\Gamma$ does
not have $\tau$ (\cite{LW}) and cannot split since groups which split
contain non abelian free groups (except of $D_{\infty}$ - the infinite
dihedral group for which  clearly $RG=0$) while amenable groups cannot
contain free groups. This last corollary was extended to all finitely
generated amenable groups in \cite{AJN}.

\chapter{Expanders and Geometry}

In this chapter we describe several ways in which expanders have appeared, somewhat unexpectedly, in geometry.  Most of these applications are for hyperbolic manifolds.  The background is given in \S 6.1.  Then in \S 6.2 we will give the first application: a proof given in [L3] using expanders and property $(\tau)$ of a conjecture of Thurston and Waldhausen
on positive virtual Betti number for arithmetic hyperbolic manifolds.  Then in \S6.3 we describe the attack of Lackenby on the ``virtual Haken conjecture" for hyperbolic 3-manifolds using expanders (or more precisely  the Lubotzky-Sarnak conjecture asserting that 3-manifolds hyperbolic groups do \emph{ not} have $(\tau)$).  While, as of now, this attack has not led to a complete solution of the virtual Haken conjecture, it led to some partial results and opened exciting new directions.  In particular, it shows connections between Heegaard genus of 3-manifolds and expanders.  This will be elaborated further in \S 6.4. We will show  there  another application of expanders to hyperbolic 3-manifolds.  Moreover, the notion of cost from dynamics will be related to 3-manifolds via expanders!

\section{Hyperbolic manifolds}

Let $M$ be an oriented $n$-dimensional hyperbolic manifold of finite volume.  Such a manifold is obtained from the Lie group $G = SO(n,1)$- the group of $(n+1)\times (n+1)$ real matrices preserving the quadratic form $X^2_1 + \ldots + X^2_n- X^2_{n+1}$ - in the following way:  Let $K = SO(n)$ sitting as a maximal compact subgroup of $G$ and $\Gamma$ a torsion-free lattice (i.e., discrete subgroup of finite covolume) in  $G$.  Then $\bbh^n= G/K$ is the $n$-dimension hyperbolic space and $M=\Gamma\setminus G/K$ is a hyperbolic manifold of finite volume.  All such manifolds are obtained like that.  Many geometric questions on such $M$ can be translated to group theoretic questions about $\Gamma$ which is actually isomorphic to the fundamental group of $M$ as $\bbh^n$ is contractible.

One of these questions is the following conjecture usually attributed to Thurston (though it probably goes back to Waldhausen):

\begin{conj}[Thurston-Waldhausen Conjecture]\label{61}
The manifold $M$ has a finite sheeted cover $M_0 \twoheadrightarrow M$ with positive $\beta_1(M_0):= \dim H_1(M_0, \bbr)$, i.e.non-trivial homology group.  Or, equivalently, in group theoretic terms:  every lattice $\Gamma$ in $SO(n, 1)$ has a finite index subgroup $\Gamma_0$ with $|\Gamma_0/[\Gamma_0, \Gamma_0]| = \infty$.
\end{conj}

The equivalence follows from two well-known facts:  every lattice is finitely generated and has a torsion free subgroup of finite index.  The commutator quotient is infinite iff there is a surjective map $\Gamma_0 \twoheadrightarrow \bbz$ and this happens iff the first real homology of $\Gamma_0\setminus G/K$ is non-trivial.

Let us mention right at the start another conjecture, due to Serre [Se] (which is now almost fully proved - see \S 6.2 below).

\begin{conj}[Serre Conjecture]\label{6.2}   If $\Gamma $ is an arithmetic lattice of $G=SO(n, 1)$, then $\Gamma$ has a negative answer to the congruence subgroup property.
\end{conj}

It is well known that the Thurston-Waldhausen conjecture  implies Serre's conjecture.  This can be seen in one of the following ways: (i) if $\Gamma $ has the congruence subgroup property then its profinite completion $\hat\Gamma$ is the same as the congruence completion.  The latter is a  product of compact $p$-adic analytic semisimple groups and as such, a finite index subgroup of it should have finite abelianization. Thus the same applies to $\hat\Gamma$ and $\Gamma$.  (ii) It is known that the congruence  subgroup property for $\Gamma$ implies super-rigidity (cf. [Se]) but if $\Gamma$ virtually maps onto $\bbz$ it does not have superigidity.

An intermediate step between these two conjectures is:
\begin{conj}[Lubotzky-Sarnak Conjecture]\label{63}
If $\Gamma$ is a lattice in $SO(n, 1)$ then $\Gamma$ does not have property $(\tau)$.
   \end{conj}
   Now, Thurston-Waldhausen conjecture $\Rightarrow $ Lubotzky-Sarnak conjecture $\Rightarrow$ Serre conjecture.
Indeed, if finite index subgroup of $\Gamma$ is mapped onto $\bbz$ then it clearly does not have $(\tau)$ as $\bbz$ does not have $(\tau)$.  Also, we mentioned in \S 2.5  that an arithmetic lattice $\Gamma$ always has $(\tau)$ with respect to congruence subgroups.  Thus, if $\Gamma$ does not have $(\tau)$, it must have also non-congruence subgroups.

This last observation was the key point in [L3] to be described in \S 6.2. But before going into these details, let us continue with another conjecture for the special case,  $n=3$, which is the most interesting case:

\begin{conj}[Virtual Haken Conjecture]\label{64} A finite volume hyperbolic 3-manifold $M$ is virtually Haken, i.e., has a finite sheeted cover which is Haken (also known as ``sufficiently large").
\end{conj}

Recall that Haken means that it contains an incompressible surface, i.e., a properly embedded orientable surface $S$ (other than $S^2$) with $\pi_1(S)$ injecting into $\pi_1 (M)$. It is known that $M$ is Haken iff $\pi_1(M)$ is either  mapped onto $\bbz$ or $\pi_1(M)$ is a free product with amalgam in a non-trivial way i.e. iff $\pi_1 (M)$ splits, in the terminology of \S 5.5.  From this it is clear that Thurston-Waldhausen conjecture for $n=3$ implies the virtual Haken conjecture.

\section{Thurston-Waldhausen conjecture for hyperbolic arithmetic manifolds}

The first use of expanders in geometry came out in the proof of the following result in [L3]:

\begin{theorem}\label{6.5} Conjecture \ref{61} is true for arithmetic lattices in $SO(n, 1)$ for $n\neq 3, 7$.  Namely, every finite volume $n$-dimensional arithmetic hyperbolic manifold has a finite sheeted cover with positive first Betti number if $n\neq 3, 7$.
\end{theorem}

The result covers also ``most" of the arithmetic lattices in $SO(3, 1)$ and $SO(7,1)$.  But these two cases are exceptional in having ``more" arithmetic lattices than what one finds in $SO(n,1)$ for other $n$'s.  The reasons are:  $SO(3, 1)$ is locally isomorphic to $SL_2(\bbc)$ and as such also has a complex structure, unlike all other $n$'s.  On the other hand $SO(7, 1)$ is a real form of $SO(8)$.  The later has Dynkin diagram of type $D_4$ and as such also has a graph automorphism of order 3 (``the triality of $D_4$") unlike the other $D_n$'s which have only  automorphisms  of order $2$.  The theory of ``Galois cohomology" which enables one to classify the arithmetic lattices in a given semisimple Lie group, shows that these anomalies give extra families of arithmetic lattices in $SO(3, 1)$ and $SO(7, 1)$ which do not exist for other $n$'s.  The method of proof of Theorem \ref{6.5} does not apply to these extra families.

The connection between Conjecture \ref{61} and expanders (or more precisely property $(\tau)$) is best explained via the following:

\begin{lem}[The Sandwich Lemma]\label{66}   Assume $G_1 \le G_2 \le G_3$ are three non-compact simple Lie groups and for each $i = 1, 2, 3$, $\Gamma_i$ is an arithmetic lattice in $G_i$ such that $\Gamma_1\le \Gamma_2 \le \Gamma_3$, $\Gamma_2 = G_2 \cap \Gamma_3$ and $\Gamma_1 = G_1 \cap \Gamma_3 (= G_1\cap \Gamma_2)$.  Then
\begin{enumerate}[(i)]
\item If $\Gamma_1$ has the Selberg property (i.e. property $(\tau)$ w.r.t. congruence subgroups - see Definition 2.14) and $\Gamma_3$ does not have $(\tau)$, then $\Ga_2$ has negative answer to the congruence subgroup problem (i.e., has non congruence subgroups).

\item{} If $\Ga_1$ has the Selberg property and $\Ga_3$ has a congruence subgroup $\Lambda$ with an infinite abelianization, then $\Ga_2$ also has such a congruence subgroup.
    \end{enumerate}

\end{lem}

Part (i) of the lemma follows immediately from Burger-Sarnak result (Theorem 2.17):  Indeed, the Selberg property of $\Gamma_1$ ``lifts up" to $\Ga_2$.  On the other hand $\Ga_2$ does not have $(\tau)$, as otherwise $\Ga_3$ would have.  Thus, $\Ga_2$ has $(\tau)$ but does not have Selberg.  In other worlds, the quotients of $\Ga_2$ modulo congruence subgroups give a family of expanders, while the family of all finite quotients does not. This gives (in a very non-constructive way!) a proof that there are non congruence subgroups in $\Ga_2$!  The proof of (ii) needs to go deeper into the actual proof of Burger-Sarnak result (see [L3] and [BS]).

Anyway, the point is that when $n\neq 3, 7$, the arithmetic lattices in $SO(n, 1)$ can be put to be $\Ga_2$ is such a sandwich: one takes $G_1 = SO(2, 1)\simeq SL_2 (\bbr)$ or $G_1= SO(3, 1)\simeq SL_2(\bbc)$ and one uses Jacquet-Langlands and Selberg results to ensure the Selberg property.  On the
other hand $G_3 = SU(n, 1)$ and one uses results of Kazhdan, Shimura and Borel-Wallach to ensure the needed properties of $\Ga_3$ - see [L3] and the references therein. These arguments show that if $\Ga$ is an arithmetic lattices in $SO(n, 1),  n\neq 3, 7$, it has a congruence subgroup which is mapped onto $\bbz$.  In particular, it does not satisfy $(\tau)$ (so Lubotzky-Sarnak Conjecture \ref{63} is also valid) and does not have the congruence subgroup property (so Serre's conjecture \ref{6.2} is also true for these cases.)

All these three conjectures are still open for the lattices of $SO(7, 1)$ coming from the triality effect of $D_4$.  The story of $n=3$ is more involved and more important.  This is our topic in the next section. We just mention in passing that the Thurston-Waldhausen conjecture has been proved for the \emph{ known} non-arithmetic lattices in $SO(n, 1)$ for $n \ge 4$ (see [L4]).  It is still widely open for others (if they exist at all ...).

\section{Hyperbolic 3-manifolds}

A few years ago, Marc Lackenby initiated a program to prove the virtual Haken conjecture for hyperbolic 3-manifolds (Conjecture \ref{64} above).  In his program, expanders (or property $\tau$) play a central role, as well as the notion of Heegaard splitting.   Let us recall the definition of the latter.

Let $M$ be a connected, closed, orientable and irreducible (i.e. any 2-dimensional sphere in $M$ bounds a 3-dimensional ball) 3-manifold.  We will be mainly interested in hyperbolic 3-manifolds, i.e., the case when $M = \Ga\backslash \bbh^3$ where $\Ga$ is a cocompact torsion free lattice in $G = PSL_2(\bbc)$.

A classical result asserts that every such $M$ can be decomposed as a union of two handle bodies $M = H_1 \cup H_2$ glued along their (isomorphic) boundaries, i.e. $H_1 \cap H_2 =\partial  H_1 =  \partial H_2$.  This decomposition \emph{(Heegaard Splitting)} is not unique.  The minimal number $g$ of handles in $H_1$ (or $H_2$- they are isomorphic) in such a decomposition is called the {\it Heegaard genus} of $M$, denoted $g(M)$.  Note that if $H_1$ has $g$ handles, its boundary is a closed surface of genus $g$ and Euler characteristic $2-2g$.

The following result of Lackenby [La3] shows a first connection between Heegaard genus and the Cheeger constant $h(M)$ of $M$.  See Definition 1.16 and Theorem 1.20  for the definition of Cheeger constant of $M$ and its connection with expanders.

\begin{theorem}\label{610} Let $M$ be a closed Riemannian 3-manifold with supremal sectional curvature $K <0$ (so $K=-1$ if $M$ is hyperbolic.)  Then
\be* h(M) \le \frac{8\pi(g(M)-1)}{|K| vol(M)}.\ee*
\end{theorem}

While this is a non-trivial result, the basic idea is simple: One proves that the Heegaard splitting (which is a topological
decomposition) can be carried out in such a way that the two parts $H_1$ and $H_2$ have approximately equal volumes - half of the volume of $M$.  Now, the area of the boundary, which is a surface of genus $g(M)$, is  given by the Gauss-Bonnet formula as a linear function of $g(M)$ and the Theorem can be deduced.

Given a Heegaard splitting of $M$ one can write down a presentation of $M$ with $g$ generators (say, pick a point on $ \partial H_1$ and take as generators the generators of $\pi_1(H_1)$ which is a free group on $g$ generators) and $g$ relations (obtained from non-trivial loops in $H_1$ which become homotopically trivial once $H_2$ is glued to $H_1$ along $\partial  H_2$).  In particular, one has:
\begin{prop}\label{67} $d(\pi_1(M)) \le g(M)$, i.e., the number of generators of the fundamental group of $M$ is bounded above by the Heegaard genus.
\end{prop}

For general 3-manifolds this can be a strict inequality, but:

\begin{conj}[Heegaard genus versus rank conjecture]\label{68} If $M$ is a compact hyperbolic 3-manifold, then $d(\pi_1 (M)) = g(M)$.
\end{conj}

Now, if $M_0$ is an $r$-sheeted cover of $M$, then the Heegaard splitting of $M$ can be lifted  to $M_0$ and one can deduce that $g(M_0) \le rg(M)$.  To ``renormalize" this, define:
\begin{definition}\label{69} Let $\Gamma = \pi_1(M), \; \call = \{ N_i\}$ a family of finite index subgroups of $\Gamma $ and $\{M_i\}$ the corresponding finite sheeted covers.  The \emph{infimal Heegaard gradient} $\chi^h_{\call}(M)$ of $M$ w.r.t. $\call $ is defined as: $\chi^h_{\call} (M) = \underset{i}{inf} \big\{\frac{2g(M_i) - 2}{[\Gamma: N_i]}\big\}$.
\end{definition}

One uses $2g(M_i)-2$, the negative of the Euler characteristic rather than $g(M_i)$, just for aesthetic reasons.  The reader may note the connection with the rank gradient defined in \S 5.5, especially in light of Proposition \ref{67}.

There are many examples of $M$ in which $\chi^h_{\call}(M) = 0$ where, say, $\call$ is a family of (all)
finite index subgroups of $\pi_1(M)$.  This happens for example if $M$  fibres over
a circle (or virtually fibres over a circle).  There are many examples of such hyperbolic 3-manifolds.
In fact, a conjecture attributed to Thurston suggests that every hyperbolic 3-manifold fibres over a
circle after passing to a suitable finite sheeted cover.  This is even stronger than Conjecture \ref{61}
above.  In group theoretic terms it is equivalent to the assertion that $\Gamma= \pi_1(M)$ has a finite index
subgroup $\Gamma_0$ with an epimorphism $\pi:\Gamma_0\twoheadrightarrow \bbz$ whose kernel $Ker(\pi)$ is finitely
generated. In this case $Ker(\pi)$ must be a surface group of genus $g$ and going along the cyclic covers
 between $\Gamma_0$ and $Ker(\pi)$ we get infinitely many covers with the same Heegaard genus.  Hence
 $\chi^h_{\call} (M) = 0$.

Lackenby conjectured that this is the only reason for the Heegaard gradient to vanish, i.e.,
\begin{conj}[Heegaard gradient conjecture]\label{611} If $M$ has a family of covers $\call$ with $\chi^h_{\call} (M) = 0$, then $M$ virtually fibres over a circle.\end{conj}

He then proved:
\begin{theorem}\label{612} Let $M$ be a closed, orientable 3-manifold and $\call = \{ N_i\} $ a family finite index normal subgroups, with corresponding covers $\{ M_i\}$.  Assume:\\
a) $\chi^h_\call(M) > 0 $ \\
b) $\Gamma = \pi_1(M)$ does \emph{not} have property $(\tau)$ w.r.t $\call$.

Then $M$ is virtually Haken.
\end{theorem}

This last result implies:
\begin{cor}\label{613} The Lubotzky-Sarnak Conjecture (Conjecture \ref{63}) and the Heegaard gradient Conjecture (Conjecture \ref{611})
imply the virtual Haken Conjecture (Conjecture \ref{64}.)
\end{cor}

Indeed let $M$ be a 3-dimensional hyperbolic manifold and $\call$ the family of all its finite index normal subgroups.  By the Lubotzky-Sarnak Conjecture $\Ga = \pi_1(M)$ does not have $(\tau)$-so condition (b) of Theorem \ref{612} is satisfied. If $\chi^h_\call (M) > 0$ then $M$ is virtually Haken by this Theorem, while if $\chi^h_\call (M) = 0$ then by the Heegaard gradient Conjecture, $M$ virtually fibres over a circle and in particular virtually Haken.

This puts property $(\tau)$ and expanders at the heart of the theory of 3-dimensional manifolds.  As of now, the full virtual Haken conjecture  has not been proven, but this approach, beside its intrinsic interest,  led to some interesting unconditional results - cf. \cite{La3}, \cite{La4}, \cite{LaLR1}, \cite{LaLR2}.

\section{Heegaard splitting, property $(\tau)$ and cost}

In the previous section we saw the result of Lackenby, Theorem \ref{610}, which connects Heegaard genus and the Cheeger constant.  The Cheeger constant is the geometric way to express expanders (see  Theorem 1.20).  These connections enabled Long, Lubotzky and Reid [LLR] to deduce the following geometric application of the theory of expanders.

\begin{theorem} \label{614} Let $M$ be a closed hyperbolic 3-manifold.  Then there exists a sequence $\call = \{ N_i\}_{i\in\bbn}$ of finite index normal subgroups of $\Gamma = \pi_1 (M)$, with $N_1\supseteq N_1\supseteq\cdots $ and\  $ \cap N_i = \{e\}$, and $\chi^h_\call (M)> 0$.  Namely, there is a constant $c>0$ such that for every $i \in\bbn$ the Heegaard genus $g(M_i) \ge c[M_i: M]$ where $M_i$ is the cover of $M$ corresponding to $N_i$ of degree $[M_i:M] = [\Gamma:N_i]$.
\end{theorem}

\begin{rem} The formulation in [LLR] is slightly weaker than what is stated here.  At that time we used the theory of sum-product results in finite fields and its applications to expanders.  The more recent results in [GSV] (see Theorem 2.24 above) enable to deduce the stronger version here.
\end{rem}

It should be stressed that in many examples of $M$'s as in Theorem \ref{614}, (and if the Thurston-Waldhausen Conjecture \ref{61} is correct, then in all such $M$'s!)  one can  also find a chain of normal subgroups $\call'$ in $\pi_1(M)$ with $\chi^h_{\call'} (M) = 0$.  This shows that the Heegaard gradient does depend on the choice of chains of normal covers (even chains with trivial intersections).

\bigskip

This brings us to another fascinating connection:  the notion of {\bf cost}.

Let $\Gamma$ be a countable group acting ergodically on $X$, a standard Borel space, by Borel automorphisms preserving a probability measure $\mu$ on $X$.  Define the equivalence relation $E$ on $X$ by $xEy$ iff $x$ and $y$ are on the same $\Gamma $-orbit.  So $E$ is a subset of $X \times X$, which can be thought as defining a graph on $X$.  For an arbitrary Borel subset $S$ of $X\times X$ we denote $deg_S(x) = | \{ y \in X| (x,y) \in S\}|$ and $e(s) = \int\limits_{x\in X} \deg_S (x) d\mu$. (See \cite{Gab2}).

We say that $S$ spans $E$ if $E$ is the minimal equivalence relation on $X$ which contains $S$ and define cost$(E) = $cost$(\Gamma, X)$ as $\inf e(S) $ when $S$ runs over all the Borel subgraphs $S$ spanning $E$.

One can easily see that if $\{ \gamma_1,\ldots,\gamma_d\}$ generates $\Gamma$ then\\ $S = \bigcup\limits^d_{i = 1} \bigcup_{x \in X} \{ (x,\gamma_i x)\}$ spans $E$ and so always cost$(\Gamma, X)\le d (\Gamma)$ - the number of generators of $\Gamma$.

This notion was introduced by Levitt [Le] and was used by Gaboriau \cite{Gab1} to distinguish between equivalence relations of different group actions.  Gaboriau conjectured:
\begin{conj}[Fixed Price Conjecture]\label{615} Given $\Gamma$, the$n$  $Cost (\Gamma, X)$ is the same number for all ergodic, essentially free actions of $\Gamma$ on a standard Borel space $X$.  (Essentially free means that the set of $x \in X$ with non-trivial stabilizer in $\Gamma$ is of measure zero).  \end{conj}

Gaboriau proved this conjecture for various groups but it is still widely open for general $\Gamma$.

An interesting example in which the cost was computed explicitly is:
\begin{theorem}[Abert-Nikolov {[AN]}]\label{616}  Let $\Gamma$ be a finitely generated group, $\call = \{ N_i\}_{i\in\bbn}$ a chain of finite index normal subgroups $N_1\supseteq N_2 \ge \cdots $ with $\underset {i}{\cap} N_i = \{ e\}$ and $\bar \Gamma = \underset{\leftarrow}{\lim} \Gamma/N_i$, the profinite completion of $\Gamma$ w.r.t. $\call$. The  group $\Gamma$ acts freely on $\bar \Gamma$ and
\be* cost(\Gamma, \bar\Gamma) = RG_\call (\Gamma) + 1\ee*
where $RG_\call (\Gamma)$ is the rank gradient of $\Gamma$ (see \S 5.5) w.r.t. $\call$, i.e., $\underset {i}{\lim} \frac{d(N_i)-1}{[\Gamma:N_i]}$.
\end{theorem}

Now, if the Fixed Price Conjecture \label{615} is true it would follow that the rank gradient of $\Ga$ does not depend on the chain $\call$ in the last theorem.  On the other  hand we saw above that for $\Ga = \pi_1 (M)$, $M$ a 3-dimensional hyperbolic compact manifold the Heegaard genus gradient \emph{ does} depend on the choice of $\call$.  This enabled Abert and Nikolov to deduce:
\begin{theorem} \label{617} At least one of the two conjectures: the Heegaard genus versus rank conjecture (Conjecture 6.9) and the Fixed price conjecture (Conjecture 6.16) is not true!
\end{theorem}

It is quite interesting how these two seemingly unrelated conjectures contradict each other  and for our story it is also interesting how this contradiction is via property $(\tau)$.

Of course, it might be that both conjectures are false! One may even speculate that the fixed price conjecture is false in general but it is true for hyperbolic groups, just as it is true for free groups (\cite{Gab1}).  If this is the case, or even if it is true for the much smaller class of fundamental groups of compact hyperbolic 3-manifolds, then the Heegaard genus versus rank conjecture would be refuted.

\chapter{Miscellaneous}

As mentioned in the introduction, expander graphs have a huge number of applications in computer science which we have not even begun to mention here.  We  have focused on applications to pure mathematics.  Even in this direction we were not able to give a comprehensive survey.   In this final chapter we will just give a list of topics that for lack of time, space or the author's expertise have not found their way into the main chapters.

\bigskip

\noindent{\bf (I) \  The Baum-Connes conjecture:}  This is a famous deep conjecture.  For a user-friendly introduction see \cite{Va3}.  Counterexamples  to a generalized form of it were given in \cite{Gro1} and \cite{HLS}. The original conjecture is still open though it was proved for many classes of groups.  The counterexamples were given by random groups constructed via expanders.

\medskip

\noindent{\bf (II) \  Embedding metric spaces:}  There is  great  interest in embedding (finite) metric spaces into Hilbert spaces in a way that the metric is more or less preserved.  In recent years this area has found many applications in computer science - see \cite[Chap. 13]{HLW}.  Expander graphs play the role of graphs whose metric is the farthest away from euclidean.

\medskip

\noindent{\bf (III) \  Dimension expanders:} The notion of expander graphs have an analogue in vector spaces.  For a fixed field $F$ and $0<\vare$, we say that $T_1,\ldots, T_k \in End_F(F^n)$, i.e. $k$ linear transformations, form an $\vare$-\emph{dimension  expander}, if for every subspace $W$ of $F^n$ with $\dim(W) \le \frac nn$, $\dim \big(\sum^k_{i=1} T_i(W)\big) \ge (1+\vare)\dim W$.  For motivation - see \cite{DS}.  Again, when one can talk on ``probability" (e.g., if $F$ is a finite or local field), ``random" $T_1,\ldots, T_k \in End_F(F^n)= M_n(F)$ will give rise to dimension expanders.  Wigderson asked for explicit constructions which are more difficult to be constructed.  This was done in \cite{LZ} for characteristic zero fields and in \cite{B1} for the general case. This motivates study of ``algebras with property $(\tau)$" like amenable algebras in \cite{Ba} and \cite{E}.

\medskip

\noindent{\bf (IV) \  High dimensional expanders:} \ A natural problem, which has been mentioned for a good number of years, is ``what is the natural definition for higher dimensional expanders?"  A suggestion for such a definition was given in \cite{Gro3} (which formally speaking does not reduce to expander for dimension one, but it still keeps  the  spirit of expanders).  In \cite{FGLNP}, random and explicit constructions of such high dimensional expanders are given.  The latter is based on \cite{LSV2}.

 \medskip

\noindent{\bf (V) \  The distribution of integer points on spheres:} \  The set of integral solutions $H_d = \{( x, y, z) \in \bbz^3 | x^2+y^2+z^2 = d\}$ can be normalized by dividing by $\sqrt{d}$ to give a subset of the sphere $S^2$.  The distribution of these points on the sphere was studied by Linnik and a modern treatment with stronger results is given in \cite{EMV}.   The modern approach makes use of random walks on expander graphs.

\medskip

\noindent{\bf (VI) \  Counting rational solutions on curves:} \  Expanders are used in a surprising way in \cite{EHK} to show some strong finiteness results on the number of $k$-rational points on some families of curves over number fields of bounded degree.

\medskip

\noindent{\bf (VII) \  $C^*$-algebras:} \   For a Hilbert space $H$, denote by $B(H)$ the $C^*$-algebra of the bounded operators of $H$. In \cite{Va1}, Ramanujan graphs were used to study the different possible norms on $B(H)\otimes B(H)$.  In \cite{BeSz}, property $(\tau)$ is used to give explicit examples of $n\times n$  matrices of norm 1 which cannot be well approximated by matrices which decompose into direct sums of smaller matrices.

In another direction, property $(\tau)$ have been used to study the question whether the set of finite dimensional representations of the $C^*$-algebra $C^*(\Gamma)$ of a finitely generated group, separate the points of $C^*(\Gamma)$ (see \cite{Be} and \cite{LSh}).

\medskip

\noindent{\bf (VIII) \  Random 3-manifolds:} \   In \cite{DT}, Dunfield  and Thurston presented a model for ``random 3-manifolds".  It is based on the fact (explained in \S 6.3) that every 3-manifold $M$ has a (non-unique) Heegaard splitting, i.e., obtained by gluing two handle-bodies along their boundaries.  The elements of the mapping class group $MCG(g)$ of a surface of genus $g$ give rise to 3-manifolds of Heegaard genus at most $g$.  Random walks on $MCG(g)$ give therefore ``random" 3-manifolds.  The group sieve method presented in Chapter 5 has already been used for studying the group $MCG(g)$  and in \cite{Ko1} and \cite{Ko2} it is used to give some results on the first homology of 3-manifolds. It  seems to have a great potential for studying further properties of ``random 3-manifolds".

We hope to return to this topic in the future.

\small
\baselineskip0pt


\begin{thebibliography}{BHKLS}


\bibitem[AJN]{AJN}  M. Abert, A. Jaikin-Zapirain and  N, Nikolov,   \emph{The rank gradient from a combinatorial viewpoint},  arXiv:math/0701925


 \bibitem[AN]{AN}  M. Abert and N. Nikolov, \emph{Rank gradient, cost of groups and the rank versus Heegaard
     genus problem}, arXiv:math/0701361


\bibitem[AC]{AC} N. Alon and  Fan R.K. Chung, \emph{Explicit construction of linear sized
tolerant networks}, Discrete Mathematics {\bf 306} (2006), 1068--1071.



\bibitem[ALW]{ALW} N. Alon, A. Lubotzky and A, Wigderson, \emph{Semi-direct product
in groups and zig-zag product in graphs: connections and applications}
(extended abstract), 42nd IEEE Symposium on Foundations of Computer
Science (Las Vegas, NV, 2001), 630637, IEEE Computer Soc., Los Alamitos, CA, 2001.





\bibitem[BHKLS]{BHKLS} L. Babai, G. Hetyei, W.M. Kantor, A. Lubotzky, and  A. Seress,  \emph{On
the diameter of finite groups},  31st Annual Symposium on Foundations of Computer Science, Vol. I, II (St. Louis, MO, 1990),  857--865, IEEE
Comput. Soc. Press, Los Alamitos, CA, 1990.

\bibitem[BKL]{BKL} L. Babai, W.M. Kantor and A. Lubotzky,
{\it  Small-diameter Cayley graphs
for finite simple groups},  European J. Combin.  10  (1989),  no. 6, 507--522.

\bibitem[BNP]{BNP} L. Babai, N. Nikolov and L. Pyber,
 \emph{Product growth and mixing
in finite groups}, Proceedings of the Nineteenth Annual ACM-SIAM Symposium
on Discrete Algorithms, 248--257, ACM, New York, 2008.


\bibitem[Ba]{Ba} L. Bartholdi,  {\it On amenability of group algebras. I},  Israel J. Math.  168  (2008), 153--165.

\bibitem[BMNVW]{BMNVW}  F. Bassino, A. Martino,  C. Nicaud,
E. Ventura and  P. Weil,
\emph{Statistical properties of subgroups of free groups},
arXiv:1001.4472

\bibitem[Be]{Be} M.B. Bekka, \emph{ On the full $C^*$-algebras of arithmetic groups and
the congruence subgroup problem}, Forum Math. 11 (1999), no. 6, 705--715.


\bibitem[BeSz]{BeSz} E.J. Benveniste  and S.J. Szarek, \emph{Property $T$, property $\tau$  and
irreducibility of matrices}, preprint.


\bibitem[B1]{B1}  J. Bourgain,  {\it Expanders and dimensional expansion}, C. R. Math. Acad.
Sci. Paris 347 (2009), no. 7-8, 357--362.

\bibitem[B2]{B2}  J. Bourgain, \emph{New developments in combinatorial number theory and
applications}, European Congress of Mathematics, 233--251, Eur. Math. Soc., Zurich, 2010.

\bibitem[BF]{BF}  J. Bourgain and  E. Fuchs,
 \emph{A proof of the positive density conjecture for integer Apollonian circle packings},
arXiv:1001.3894


\bibitem[BFLM]{BFLM}  J. Bourgain, A. Furman, E. Lindenstrauss and S. Mozes,  \emph{Invariant measures and stiffness for non-abelian groups of toral automorphisms}, C. R. Math. Acad. Sci. Paris 344 (2007), no. 12, 737--742.

\bibitem[BG1]{BG1} J. Bourgain and A. Gamburd, {\it Uniform expansion bounds for Cayley
graphs of ${\rm SL}_2(\Bbb F_p)$},  Ann. of Math. (2)  167  (2008),  no. 2, 625--642.

\bibitem[BG2]{BG2} J. Bourgain and A. Gamburd, \emph{On the spectral gap for
finitely-generated subgroups of $\rm SU(2)$}, Invent. Math. 171 (2008), no. 1, 83--121.

\bibitem[BG3]{BG3}  J. Bourgain and A. Gamburd, \emph{Expansion and random walks in ${\rm
SL}_d(\Bbb Z/p^n\Bbb Z)$. I}, J. Eur. Math. Soc. (JEMS) 10 (2008), no. 4, 987--1011.

\bibitem[BG4]{BG4} J. Bourgain and A. Gamburd, \emph{Expansion and random walks in ${\rm
SL}_d(\Bbb Z/p^n\Bbb Z)$. II}, with an appendix by J. Bourgain,
J. Eur. Math. Soc. (JEMS) 11 (2009), no. 5, 1057--1103.


\bibitem[BGS1]{BGS1}  J. Bourgain, A. Gamburd and P. Sarnak, \emph{Sieving and
expanders}, C. R. Math. Acad. Sci. Paris 343 (2006), no. 3, 155--159.

\bibitem[BGS2]{BGS2} J. Bourgain, A. Gamburd and P. Sarnak, {\it Affine linear sieve,
expanders, and sum-product}, Invent. Math. 179 (2010), no. 3, 559--644.

\bibitem[BGS3]{BGS3} J. Bourgain, A. Gamburd and P. Sarnak,
{\it Generalization of Selberg's 3/16 Theorem and Affine Sieve}, arXiv:0912.5021

\bibitem[BKT]{BKT} J. Bourgain, N. Katz and T. Tao, \emph{A sum-product estimate in finite
fields, and applications}, Geom. Funct. Anal. 14 (2004), no. 1, 27--57.

\bibitem[BK]{BK} J. Bourgain and  A. Kontorovich,
\emph{On representations of integers in thin subgroups of $SL(2,Z)$}, Geom.~Funct.~Anal.~(GAFA) {\bf 20} (2010), 1144--1174.

\bibitem[BV]{BV} J. Bourgain and P.P. Varju,
{\it Expansion in $SL_d(Z/qZ)$, $q$ arbitrary},
arXiv:1006.3365

\bibitem[BCLM]{BCLM} E. Breuillard, Y. De Cornulier, A. Lubotzky and C. Meiri,  \emph{Conjugacy
growth of linear groups}, preprint.

\bibitem[BGa]{BGa} E. Breuillard and A. Gamburd,
{\it Strong uniform expansion in $\mathrm{SL}(2,p)$},
 Geom.\ Funct.\ Anal.\ (GAFA) {\bf 20} (2010), 1201--1209.

\bibitem[BGT1]{BGT1}
E. Breuillard, B. Green and  T. Tao,
{\it Linear Approximate Groups},
 Electron.~Res.~Announc.~Math.~Sci. {\bf 17} (2010), 57--67.

\bibitem[BGT2]{BGT2}  E. Breuillard, B. Green and T. Tao,
{\it Approximate subgroups of linear groups},
 arXiv:1005.1881

\bibitem[BGT3]{BGT3}  E. Breuillard, B. Green and T. Tao,
{\it Suzuki groups as expanders},
arXiv:1005.0782

\bibitem[BGGT1]{BGGT1}  E. Breuillard, B. Green, R. Guralnick and  T. Tao,
\emph{Strongly dense free subgroups of semisimple algebraic groups},
arXiv:1010.4259

\bibitem[BGGT2]{BGGT2} E. Breuillard, B. J. Green, R. Guralnick and T. C. Tao, \emph{Expansion in
finite simple groups of Lie type}, in preparation.

\bibitem[BLMS]{BLMS}  Y. Bugeaud, F. Luca,  M. Mignotte and S. Siksek, \emph{ On
Fibonacci numbers with few prime divisors}, Proc. Japan Acad. Ser. A Math. Sci. 81 (2005), no. 2, 17--20.

\bibitem[BS]{BS}  M. Burger and P. Sarnak, \emph{Ramanujan duals. II}, Invent. Math. 106
(1991), no. 1, 111.

\bibitem[Bu]{Bu} Peter  Buser, \emph{ Geometry and spectra of compact Riemann surfaces},
Progress in Mathematics, 106. Birkhuser Boston, Inc., Boston, MA, 1992. xiv+454 pp.

\bibitem[CLMNO]{CLMNO} F. Celler, C.R. Leedham-Green, S.H. Murray,
A.C. Niemeyer and E.A. O'Brien, \emph{Generating random elements of a finite
group}, Comm. Algebra 23 (1995), no. 13, 4931--4948.

\bibitem[Ch]{Ch} J.R. Chen,   \emph{On the representation of a larger even integer as the
sum of a prime and the product of at most two primes}, Sci. Sinica 16
(1973), 157--176.

\bibitem[Cl]{Cl} L. Clozel,  \emph{D\'emonstration de la conjecture $\tau$}, Invent. Math.
151 (2003), no. 2, 297--328.

\bibitem[DSV]{DSV}  G. Davidoff, P. Sarnak and A. Valette, \emph{Elementary number
theory, group theory, and Ramanujan graphs}, London Mathematical Society
Student Texts, 55. Cambridge University Press, Cambridge, 2003. x+144 pp.

\bibitem[DSC]{DSC} P. Diaconis and L. Saloff-Coste, \emph{ Walks on generating sets of groups},
Invent. Math. 134 (1998), no. 2, 251--299.

\bibitem[Di]{Di}   O. Dinai, \emph{Expansion properties of finite simple groups}, arXiv:1001.5069

\bibitem[D]{D} J.D. Dixon, {\it The probability of generating the symmetric group},
Math. Z.  110  1969 199--205.

\bibitem[DT]{DT}  N.M. Dunfield and W.P. Thurston, {\it Finite covers of random
3-manifolds}, Invent. Math. 166 (2006), no. 3, 457--521.

\bibitem[DS]{DS} Z. Dvir, and A. Shpilka, \emph{Towards dimension expanders over finite
fields}, Twenty-Third Annual IEEE Conference on Computational Complexity,
304-310, IEEE Computer Soc., Los Alamitos, CA, 2008.

\bibitem[E]{E} G. Elek, {\it The amenability of affine algebras},  J. Algebra  264 (2003),  no. 2, 469--478.

\bibitem[EHK]{EHK}  J. Ellenberg, C. Hall and E. Kowalski,
{\it Expander graphs, gonality and variation of Galois representations},
 arXiv:1008.3675

\bibitem[EMV]{EMV}  J. S. Ellenberg, P. Michel and A. Venkatesh, \emph{Linnik's ergodic method and the distribution of integer points
on spheres},
 arXiv:1001.0897

\bibitem[Er]{Er} M. Ershov, \emph{Golod-Shafarevich groups with property $(T)$ and
Kac-Moody groups}, Duke Math. J. 145 (2008), no. 2, 309--339.

\bibitem[EJ]{EJ} M. Ershov and A. Jaikin-Zapirain, \emph{ Property $(T)$ for
noncommutative universal lattices}, Invent. Math. 179 (2010), no. 2, 303--347.

\bibitem[FGLNP]{FGLNP}  J. Fox, M. Gromov, V. Lafforgue, A. Naor and  J. Pach,
    \emph{Overlap properties of geometric expanders},
    arXiv:1005.1392

\bibitem[FI]{FI} J. Friedlander and H. Iwaniec, \emph{ Opera de cribro},  American
Mathematical Society Colloquium Publications, 57. American Mathematical
Society, Providence, RI, 2010. xx+527 pp.

\bibitem[Fr]{Fr}  J. Friedman, \emph{ A proof of Alon's second eigenvalue conjecture and
related problems}, Mem. Amer. Math. Soc. 195 (2008), no. 910, viii+100 pp

\bibitem[Fu]{Fu} E. Fuchs, Ph.D Thesis, Princeton University.

\bibitem[FS]{FS}  E. Fuchs and K. Sanden,
   \emph{ Some experiments with integral Apollonian circle packings},
    arXiv:1001.1406

\bibitem[Gab1]{Gab1} D. Gaboriau, \emph{ Co\^{u}t des relations d'\'equivalence et des groupes},
 Invent. Math. 139 (2000), no. 1, 41--98.

\bibitem[Gab2]{Gab2} D. Gaboriau, \emph{What is Cost?} Notices AMS {\bf 57} (2010), 1295--1296.

\bibitem[Ga1]{Ga1}  A. Gamburd, \emph{On the spectral gap for infinite index ``congruence"
subgroups of ${\rm SL}_2(\bold Z)$}, Israel J. Math. 127 (2002), 157--200.

\bibitem[Ga2]{Ga2}  A. Gamburd, \emph{Expander graphs, random matrices and quantum chaos},
Random walks and geometry, 109--140, Walter de Gruyter GmbH \& Co. KG, Berlin, 2004.

\bibitem[GHSSV]{GHSSV} A.  Gamburd, S. Hoory, M. Shahshahani. A. Shalev and B. Virg, \emph{ On
the girth of random Cayley graphs}, Random Structures Algorithms 35 (2009),
no. 1, 100--117.

\bibitem[GJS]{GJS}  A. Gamburd, D. Jakobson and P. Sarnak, \emph{Spectra of elements
in the group ring of ${\rm SU}(2)$}, J. Eur. Math. Soc. (JEMS) 1 (1999), no. 1, 51--85.

\bibitem[GH]{GH} N. Gill and H.A. Helfgott,
\emph{ Growth of small generating sets in $SL_n(Z/pZ)$},
arXiv:1002.1605

\bibitem[Gl]{Gl} G. Glauberman,  \emph{Factorizations in local subgroups of finite groups},
Regional Conference Series in Mathematics, No. 33. American Mathematical Society, Providence, R.I., 1977. ix+74 pp

\bibitem[Go]{Go} W.T. Gowers, \emph{ Quasirandom groups}, Combin. Probab. Comput. 17 (2008), no. 3, 363--387.

\bibitem[GLMWY]{GLMWY}  R.L. Graham, J.C. Lagarias, C.L. Mallows, L. Colin, A.R. Wilks and
C.H.  Yan, \emph{Apollonian circle packings: number theory}, J. Number Theory 100 (2003), no. 1, 1--45.

\bibitem[Gr]{Gr} B. Green,
     \emph{Approximate groups and their applications: work of Bourgain,
Gamburd, Helfgott and Sarnak},
    arXiv:0911.3354

\bibitem[GT1]{GT1} B. Green and T. Tao, \emph{ The primes contain arbitrarily long
arithmetic progressions}, Ann. of Math. (2) 167 (2008), no. 2, 481--547.

\bibitem[GT2]{GT2} B. Green and T. Tao, \emph{ Linear equations in primes}, Ann. of
Math. (2) 171 (2010), no. 3, 1753--1850.

\bibitem[GTZ]{GTZ}  B. Green, T. Tao and  T. Ziegler,
     \emph{An inverse theorem for the Gowers $U^{s+1}[N]$-norm},
     arXiv:1009.3998



\bibitem[Gro1]{Gro1} M. Gromov, \emph{ Random walk in random groups}, Geom. Funct. Anal. 13
(2003), no. 1, 73--146.

\bibitem[Gro2]{Gro2} M. Gromov, \emph{Singularities, expanders and topology of maps. I.
Homology versus volume in the spaces of cycles}, Geom. Funct. Anal. 19
(2009), no. 3, 743--841.

\bibitem[Gro3]{Gro3} M. Gromov, \emph{Singularities, expanders and topology of maps. II.
From combinatorics to topology via
algebraic isoperimetry}, Geom. Funct. Anal.\ {\bf 20} (2010), 416--526.

\bibitem[GrGu]{GrGu}  M. Gromov and  L. Guth,
\emph{Generalizations of the Kolmogorov-Barzdin embedding estimates},
    arXiv:1103.3423

\bibitem[GL]{GL} F. Grunewald and A. Lubotzky, {\it Linear representations of the
automorphism group of a free group}, Geom. Funct. Anal. 18 (2009), no. 5, 1564--1608.

\bibitem[HL]{HL} G.H. Hardy and J.E. Littlewood, \emph{ Some problems of `Partitio   numerorum'; III: On the expression of a number as a sum of primes}, Acta
Math. 44 (1923), no. 1, 1--70.

\bibitem[H]{H}  H.A. Helfgott, {\it Growth and generation in ${\rm SL}_2(\Bbb Z/p\Bbb
Z)$}, Ann. of Math. (2) 167 (2008), no. 2, 601--623.

\bibitem[HLS]{HLS}   N. Higson,  V. Lafforgue and G. Skandalis, \emph{Counterexamples to the
Baum-Connes conjecture}, Geom. Funct. Anal. 12 (2002), no. 2, 330--354.

\bibitem[HLW]{HLW}  S. Hoory, N. Linial and A. Wigderson, {\it     Expander graphs and
their applications}, Bull. Amer. Math. Soc. (N.S.) 43 (2006), no. 4, 439--561.

\bibitem[HKLS]{HKLS}  E. Hrushovski, P.H. Kropholler, A. Lubotzky and A. Shalev, {\it  Powers
in finitely generated groups}, Trans. Amer. Math. Soc. 348 (1996), no. 1, 291--304.

\bibitem[IK]{IK}  H. Iwaniec and E. Kowalski, \emph{ Analytic number theory}, American
Mathematical Society Colloquium Publications, 53. American Mathematical
Society, Providence, RI, 2004. xii+615 pp.

\bibitem[JKZ]{JKZ}   F. Jouve, E. Kowalski and  D. Zywina,
{\it Splitting fields of characteristic polynomials of random elements
in arithmetic groups}, Israel J.~of Math., to appear.
arXiv:1008.3662

\bibitem[KL]{KL} W.M. Kantor and A. Lubotzky, {\it  The probability of generating
a finite classical group},  Geom. Dedicata  36  (1990),  no. 1, 67--87.

\bibitem[KLS]{KLS} W. M. Kantor, A. Lubotzky and A. Shalev,
    \emph{Invariable generation and the Chebotarev invariant of a finite group},
    arXiv:1010.5722

\bibitem[KMSS1]{KMSS1} I. Kapovich, A. Miasnikov, P. Schupp and V. Shpilrain,
{\it  Generic-case complexity, decision problems in group theory, and
random walks}, J. Algebra 264 (2003), no. 2, 665--694.

 \bibitem[KMSS2]{KMSS2}  I. Kapovich, A. Miasnikov, P. Schupp and V. Shpilrain, {\it Average-case complexity and decision problems in group theory},
Adv. Math. 190 (2005), no. 2, 343--359.

\bibitem[KS1]{KS1} I. Kapovich and P. Schupp, \emph{ On group-theoretic models of randomness
and genericity}, Groups Geom. Dyn. 2 (2008), no. 3, 383--404.

\bibitem[K1]{K1} M. Kassabov, {\it Universal lattices and unbounded rank expanders},
Invent. Math. 170 (2007), no. 2, 297--326.

\bibitem[K2]{K2} M. Kassabov,  {\it Symmetric groups and expander graphs}, Invent. Math.
170 (2007), no. 2, 327--354.


\bibitem[KLN]{KLN} M. Kassabov, A. Lubotzky and N. Nikolov,  {\it Finite simple groups as expanders}, Proc. Natl. Acad. Sci. USA 103 (2006), no. 16,
6116--6119.

\bibitem[KN]{KN} M. Kassabov and N.  Nikolov, {\it Universal lattices and property
tau}, Invent. Math. 165 (2006), no. 1, 209--224.

\bibitem[KaL]{KaL}  T. Kaufman and  A. Lubotzky,
\emph{Edge transitive Ramanujan graphs and
    highly  symmetric LDPC good codes}, preprint.

\bibitem[KaW]{KaW} T. Kaufman and A. Wigderson,
 \emph{Symmetric LDPC and local Testing},
     Innovations in Computer Science, 406--421, 2010.


\bibitem[Ka]{Ka} D.A. Kazhdan, \emph{On the connection of the dual space of a group with
the structure of its closed subgroups}, (Russian) Funkcional. Anal. i Priloen. 1 1967,  71--74.

\bibitem[Ki]{Ki} H.H. Kim,  \emph{Functoriality for the exterior square of ${\rm GL}_4$
and the symmetric fourth of ${\rm GL}_2$}, with appendix 1 by Dinakar
Ramakrishnan and appendix 2 by Kim and Peter Sarnak. J. Amer. Math. Soc.
16 (2003), no. 1, 139--183.

\bibitem[KB]{KB} A. N.Kolmogorov and Y.M. Barzdin, \emph{On the realization of nets in
3-dimensional space}, Probl. Cybernet, 8, 261--268, 1967. See also
Selected Works of A.N. Kolmogorov, Vol 3, pp 194--202 (and a remark on page
245), Kluwer Academic Publishers, 1993.

\bibitem[KO1]{KO1}  A. Kontorovich  and H. Oh,
    \emph{Apollonian circle packings and closed horospheres on hyperbolic 3-manifolds},
    arXiv:0811.2236

\bibitem[KO2]{KO2}  A. Kontorovich  and H. Oh,
\emph{ Almost prime Pythagorean triples in thin orbits},
arXiv:1001.0370

\bibitem[Ko1]{Ko1} E. Kowalski, {\it The large sieve and its applications}, Arithmetic geometry, random walks and discrete groups. Cambridge Tracts in
Mathematics, 175. Cambridge University Press, Cambridge, 2008. xxii+293 pp.


\bibitem[Ko2]{Ko2} E. Kowalski, \emph{Sieve and expansion}, Seminar Bourbaki, November 2010.

\bibitem[La1]{La1} M. Lackenby, \emph{Expanders, rank and graphs of groups},  Israel J. Math. 146 (2005), 357--370.

\bibitem[La2]{La2} M. Lackenby,  \emph{A characterisation of large finitely presented
groups}, J. Algebra 287 (2005), no. 2, 458--473.

\bibitem[La3]{La3} M. Lackenby, {\it Heegaard splittings, the virtually Haken conjecture
and property $(\tau)$},  Invent. Math.  164  (2006),  no. 2, 317--359.

\bibitem[La4]{La4} M. Lackenby, \emph{Large groups, property $(\tau)$ and the homology
growth of subgroups}, Math. Proc. Cambridge Philos. Soc. 146 (2009), no. 3,
625--648.

\bibitem[LaLR1]{LaLR1} M. Lackenby, D.D.  Long and A.W.  Reid,  {\it    LERF and the
Lubotzky-Sarnak conjecture},  Geom. Topol.  12  (2008),  no. 4, 2047--2056.

\bibitem[LaLR2]{LaLR2} M. Lackenby, D.D. Long and A.W.  Reid,  \emph{ Covering spaces of
arithmetic 3-orbifolds}, Int. Math. Res. Not. IMRN 2008, no. 12, Art. ID rnn036, 38 pp.

\bibitem[Le]{Le} G. Levitt, \emph{On the cost of generating an equivalence relation}, Ergodic
Theory Dynam, Systems, 15(6) (1995), 1173--1181.

\bibitem[LiSh]{LiSh} M.W. Liebeck and A. Shalev, {\it The probability of generating a
finite simple group},  Geom. Dedicata  56  (1995),  no. 1, 103--113.

\bibitem[Li]{Li} N. Linial, \emph{Finite metric-spaces: combinatorics, geometry and
algorithms}, Proceedings of the International Congress of Mathematicians,
Vol. III (Beijing, 2002), 573--586, Higher Ed. Press, Beijing, 2002.

\bibitem[LLR]{LLR} D.D. Long, A. Lubotzky and A.W.  Reid, {\it     Heegaard genus and property
$\tau$ for hyperbolic 3-manifolds}, J. Topol. 1 (2008), no. 1, 152--158.

\bibitem[Lo]{Lo} E. Looijenga, \emph{Prym representations of mapping class groups},  Geom.
Dedicata 64 (1997), no. 1, 69--83.

\bibitem[L1]{L1} A. Lubotzky,  {\it Discrete groups, expanding graphs and invariant
measures}, with an appendix by Jonathan D. Rogawski. Reprint of the 1994 edition. Modern Birkh\"auser Classics. Birkh\"auser Verlag, Basel, 2010. iii+192
pp.

\bibitem[L2]{L2} A. Lubotzky, {\it Cayley graphs: eigenvalues, expanders and random
walks}, Surveys in combinatorics, 1995 (Stirling), 155--189, London Math.
Soc. Lecture Note Ser., 218, Cambridge Univ. Press, Cambridge, 1995.

\bibitem[L3]{L3} A. Lubotzky,  {\it Eigenvalues of the Laplacian, the first Betti
number and the congruence subgroup problem}, Ann. of Math. (2) 144 (1996), no. 2, 441--452.

\bibitem[L4]{L4} A. Lubotzky, \emph{ Free quotients and the first Betti number of some
hyperbolic manifolds}, Transform. Groups 1 (1996), no. 1-2, 71--82.

\bibitem[L5]{L5}   A. Lubotzky, \emph{ What is$\dots$property $(\tau)$}, Notices Amer.
Math. Soc. 52 (2005), no. 6, 626--627.

\bibitem[L6]{L6} A. Lubotzky,
{\it Finite simple groups of Lie type as expanders},
 to appear in J. Eur. Math. Soc.
 arXiv:0904.3411

\bibitem[LM1]{LM1} A. Lubotzky and C. Meiri, \emph{  Sieve methods in group theory: I. powers in
linear groups}, preprint.

\bibitem[LM2]{LM2} A. Lubotzky and C. Meiri,  \emph{   Sieve methods in group theory: II. The
mapping class group}, preprint.

\bibitem[LP]{LP} A. Lubotzky and I. Pak, \emph{ The product replacement algorithm and
Kazhdan's property $(T)$}, J. Amer. Math. Soc. 14 (2001), no. 2, 347--363.

\bibitem[LPS1]{LPS1} A. Lubotzky, R. Phillips and P. Sarnak, \emph{ Ramanujan conjecture and
explicit construction of
expanders}, Proc. STOC. 86 (1986), 240--246.

\bibitem[LPS2]{LPS2}  A. Lubotzky, R. Phillips and P. Sarnak,  \emph{ Ramanujan graphs},
Combinatorica 8 (1988), no. 3, 261--277.

\bibitem[LPS3]{LPS3} A. Lubotzky, R. Phillips and P. Sarnak,  \emph{ Hecke operators and
distributing points on the sphere. I}, Frontiers of the mathematical sciences: 1985 (New York, 1985). Comm. Pure Appl. Math. 39 (1986), no. S,
suppl., 149--186.

\bibitem[LPS4]{LPS4} A. Lubotzky, R. Phillips and P. Sarnak, \emph{ Hecke operators and
distributing points on $S2$. II}, Comm. Pure Appl. Math. 40 (1987), no.
4, 401--420.

\bibitem[LR]{LR} A. Lubotzky and L. Rosenzweig,
 \emph{The galois groups of random elements
    of linear groups}, in preparation.

\bibitem[LSV1]{LSV1} A. Lubotzky, B.  Samuels and U. Vishne,  \emph{ Ramanujan complexes
of type $\tilde A_d$}, Probability in Mathematics. Israel J. Math. 149
(2005), 267--299.

\bibitem[LSV2]{LSV2} A. Lubotzky, B. Samuels and U. Vishne, \emph{ Explicit
constructions of Ramanujan complexes of type $\tilde A_d$}, European J.
Combin. 26 (2005), no. 6, 965--993.


\bibitem[LS]{LS} A. Lubotzky and D. Segal, \emph{Subgroup growth}, Progress in
Mathematics, 212. Birkhuser Verlag, Basel, 2003. xxii+453 pp.

\bibitem[LSh]{LSh}  A. Lubotzky and Y. Shalom, \emph{ Finite representations in the
unitary dual and Ramanujan groups}, Discrete geometric analysis, 173--189,
Contemp. Math., 347, Amer. Math. Soc., Providence, RI, 2004.

\bibitem[LW]{LW}  A. Lubotzky and B. Weiss, \emph{ Groups and expanders}, Expanding graphs
(Princeton, NJ, 1992), 95109, DIMACS Ser. Discrete Math. Theoret. Comput.
Sci., 10, Amer. Math. Soc., Providence, RI, 1993.

\bibitem[LZ]{LZ} A. Lubotzky and E. Zelmanov, {\it Dimension expanders},  J. Algebra
319 (2008), no. 2, 730--738.

\bibitem[LZi]{LZi}  A. Lubotzky and R.J.  Zimmer, \emph{ Variants of Kazhdan's
property for subgroups of semisimple groups}, Israel J. Math. 66 (1989), no. 1-3, 289--299.

\bibitem[LZu]{LZu} A. Lubtozky and A. Zuk, \emph{ On property $(\tau)$}, monograph in preparation.

\bibitem[Ma1]{Ma1} J. Maher,  {\it Random walks on the mapping class group},  arXiv:math/0604433

\bibitem[Ma2]{Ma2} J. Maher, {\it Random Heegaard splittings}, J.\ Topology {\bf 3} (2010), 997--1025.



 \bibitem[MS]{MS}  J. Malestein and  J. Souto,
\emph{On genericity of pseudo-Anosovs in the Torelli group},
    arXiv:1102.0601



\bibitem[M1]{M1} G.A. Margulis, \emph{Explicit constructions of expanders}. (Russian)
Problemy Peredac(i Informacii 9 (1973), no. 4, 7180. English translation:
Problems of Information Transmission 9 (1973), no. 4, 325--332 (1975).

\bibitem[M2]{M2} G.A. Margulis, \emph{ Explicit constructions of graphs without short cycles
and low density codes}, Combinatorica 2 (1982), no. 1, 71--78.

\bibitem[M3]{M3} G.A. Margulis, \emph{ Explicit group-theoretic constructions of combinatorial schemes and
their applications in the construction of expanders and concentrators},
Problems of
Information Transmission, 24(1):39--46, 1988.

\bibitem[MVW]{MVW} C.R.  Matthews, L.N.  Vaserstein and B. Weisfeiler, \emph{ Congruence
properties of Zariski-dense subgroups. I}, Proc. London Math. Soc. (3) 48
(1984), no. 3, 514--532.

\bibitem[Maz]{Maz} B. Mazur, \emph{It is a story}, A lecture given at Diaconis' 60th
birthday. Available at
http://www.math.ucsd.edu/~williams/diaconis/It.is.a.story.3.pdf

\bibitem[MW]{MW} R. Meshulam and A. Wigderson, \emph{Expanders in group algebras},
Combinatorica 24 (2004), no. 4, 659--680.

\bibitem[Mo]{Mo} M. Morgenstern,
\emph{Existence and explicit constructions of
    $q+1$ regular Ramanujan graphs for every prime power $q$}, J.~Combin.~Theory Ser.~B {\bf 62} (1994), 44--62.

\bibitem[MR]{MR}  A.G. Myasnikov and A.N. Rybalov, {\it Generic complexity of
undecidable problems}, J. Symbolic Logic 73 (2008), no. 2, 656--673.

\bibitem[NS]{NS}   A. Nevo and P. Sarnak,
    \emph{Prime and almost prime integral points on principal homogeneous
spaces}, Acta Math.\ {\bf 205} (2010), 361--402.

\bibitem[Ni]{Ni}  N. Nikolov, \emph{ A product of decomposition for the classical
quasisimple groups}, J. Group Theory 10 (2007), no. 1, 43--53.

\bibitem[NiPy]{NiPy}  N. Nikolov and L. Pyber,
    \emph{Product decompositions of quasirandom groups and a Jordan type theorem},
    arXiv:math/0703343



\bibitem[No]{No} M.V. Nori, \emph{ On subgroups of ${\rm GL}_n({\bf F}_p)$},
Invent. Math. 88 (1987), no. 2, 257--275.


\bibitem[Pi]{Pi} R. Pink, \emph{ Strong approximation for Zariski dense subgroups over
arbitrary global fields}, Comment. Math. Helv. 75 (2000), no. 4, 608--643.

\bibitem[Pin]{Pin} M.S. Pinsker, \emph{On the complexity of a concentrator}, 7th International Teletraffic Conference, Stockholm, pages 318/1--318/4, June
1973.

\bibitem[PS1]{PS1} L. Pyber and E. Szab\'o,
{\it Growth in finite simple groups of Lie type},
arXiv:1001.4556

\bibitem[PS2]{PS2} L. Pyber and E. Szab\'o,
{\it  Growth in finite simple groups of Lie type of bounded rank},
arXiv:1005.1858

\bibitem[RVW]{RVW} O. Reingold, S. Vadhan and A. Wigderson, \emph{ Entropy waves, the
zig-zag graph product, and new constant-degree expanders}, Ann. of Math. (2) 155 (2002), no. 1, 157--187.

\bibitem[Ri]{Ri} I. Rivin, {\it Walks on groups, counting reducible matrices,
polynomials, and surface and free group automorphisms},  Duke Math. J.  142
(2008),  no. 2, 353--379.

\bibitem[RSW]{RSW}  E. Rozenman, A. Shalev and A. Wigderson,  \emph{ Iterative construction
of Cayley expander graphs}, Theory Comput. 2 (2006), 91--120.

\bibitem[SGS]{SGS} A. Salehi-Golsefiday and P. Sarnak, \emph{Affine linear sieve}, in
preparation.

\bibitem[SGV]{SGV}  A. Salehi-Golsefiday and  P. Varju, \emph{Expansion in perfect groups},
preprint.

\bibitem[S1]{S1}  P. Sarnak, \emph{ Some applications of modular forms}, Cambridge Tracts
in Mathematics, 99. Cambridge University Press, Cambridge, 1990. x+111 pp.

\bibitem[S2]{S2} P. Sarnak,  \emph{Selberg's eigenvalue conjecture},  Notices Amer. Math.
Soc. 42 (1995), no. 11, 1272--1277.

\bibitem[S3]{S3}  P. Sarnak, \emph{ What is $\dots$ an expander}? Notices Amer. Math. Soc. 51
(2004), no. 7, 762--763.

\bibitem[S4]{S4} P. Sarnak, \emph{ Equidistribution and primes}, G\'eom\'etrie differentielle,
physique math\'ematiques, math\'ematiques et soci\'et\'e. II. Asterisque No. 322 (2008),
225--240.

\bibitem[S5]{S5} P.  Sarnak, \emph{ Letter to Lagarias on integral Apollonian
packings}. Available at
http://www.math.princeton

\bibitem[S6]{S6}  P. Sarnak, \emph{Equidistribution and Primes}, (2007) PIMS Lecture. Available at
http://www.math.princeton

\bibitem[S7]{S7} P. Sarnak, \emph{Primes and orbits}, MAA Garden State lecture. Available at
http://www.math.princeton

\bibitem[S8]{S8} P. Sarnak, \emph{Integral Apollonian Packings} - MAA Lecture January 2009.
 Available at
 http://www.math.princeton

\bibitem[SS]{SS} A. Schinzel and W. Sierpin'ski, \emph{ Sur certaines hypoth\'eses concernant les
nombres premiers}, (French) Acta Arith. 4 (1958), 185--208; erratum 5 (1958)
259.

\bibitem[SiSp]{SiSp}  M. Sipser and D.A. Spielman,
\emph{Expander codes}, IEEE Trans.\ Inform.\ Theory {\bf 42} (1996), 1710--1722.


\bibitem[Sel]{Sel} A. Selberg, \emph{On the estimation of Fourier coefficients of modular
forms}, 1965 Proc. Sympos. Pure Math., Vol. VIII pp. 115 Amer. Math. Soc., Providence, R.I.

\bibitem[Se]{Se} J-P. Serre, \emph{ Le probl\`eme des groupes de congruence pour $SL_2$},
(French) Ann. of Math. (2) 92 1970 489--527.

\bibitem[Sh1]{Sh1}  Y. Shalom, \emph{ Expanding graphs and invariant means}, Combinatorica
17 (1997), no. 4, 555--575.

\bibitem[Sh2]{Sh2} Y. Shalom, \emph{ Expander graphs and amenable quotients}, Emerging
applications of number theory (Minneapolis, MN, 1996), 571--581, IMA Vol.
Math. Appl., 109, Springer, New York, 1999.

\bibitem[Sh3]{Sh3} Y. Shalom, \emph{ Bounded generation and Kazhdan's property $(T)$}, Inst.
Hautes \'Etudes Sci. Publ. Math. No. 90 (1999), 145--168 (2001).

\bibitem[Sh4]{Sh4} Y. Shalom, \emph{ The algebraization of Kazhdan's property $(T)$},
International Congress of Mathematicians. Vol. II, 1283--1310, Eur. Math.
Soc., Zurich, 2006.

\bibitem[Ta]{Ta} R. M. Tanner,
\emph{A recursive approach to low complexity codes},
    IEEE Transactions on Information Theory {\bf 27} (1981), 533--547.

\bibitem[TV]{TV}   T. Tao and V.  Vu, \emph{Additive combinatorics}, Cambridge Studies in
Advanced Mathematics, 105. Cambridge University Press, Cambridge, 2006.
xviii+512 pp.

\bibitem[Ti]{Ti} J. Tits,  \emph{Free subgroups in linear groups}, J. Algebra 20 1972 250--270.

\bibitem[Va1]{Va1} A. Valette, \emph{An application of Ramanujan graphs to $C^*$-algebra
tensor products}, 15th British Combinatorial Conference (Stirling, 1995).
Discrete Math. 167/168 (1997), 597--603.

\bibitem[Va2]{Va2} A.  Valette, \emph{ Graphes de Ramanujan et applications}, (French)
\emph{Ramanujan graphs and applications}, S\'eminaire Bourbaki, Vol. 1996/97.
Astrisque No. 245 (1997), Exp. No. 829, 4, 247--276.

\bibitem[Va3]{Va3} A. Valette, \emph{ Introduction to the Baum-Connes conjecture},   lectures
in Mathematics ETH Zrich. Birkhuser Verlag, Basel, 2002.

\bibitem[V]{V}  P.O. Varju, {\it Expansion in $SL_d(O_K/I)$, $I$ square-free},
arXiv:1001.3664.

\bibitem[W]{W} B.  Weisfeiler, \emph{Strong approximation for Zariski-dense subgroups of
semisimple algebraic groups}, Ann. of Math. (2) 120 (1984), no. 2, 271--315.

\bibitem[Z]{Z} A. Zuk, \emph{ Property (T) and Kazhdan constants for discrete groups}, Geom.
Funct. Anal. 13 (2003), no. 3, 643--670.

\end{thebibliography}
\end{document}